\documentclass[english,11pt,letterpaper]{article}
\usepackage[T1]{fontenc}
\usepackage{babel}
\usepackage{csquotes}
\usepackage{graphicx}
\usepackage{caption}
\usepackage{subcaption}

\usepackage{lscape}
\usepackage{todonotes}
\usepackage[margin=1in]{geometry}
\usepackage{amsmath}
\usepackage{amssymb}
\usepackage{mathtools}
\usepackage{booktabs}
\usepackage{enumitem}
\usepackage{longtable}
\usepackage{makecell}

\newcommand{\var}[1]{{\color{red} #1}}

\usepackage[backend=biber,citestyle=numeric-comp,bibstyle=ieee,sorting=nyt,minbibnames=5,maxbibnames=5,defernumbers=true,giveninits,doi=false,isbn=false,url=false,eprint=false]{biblatex}
\title{A unified vertical alignment and earthwork model in road design with a new convex optimization model for road networks\thanks{This is an original manuscript of an article published by Taylor \& Francis in Engineering Optimization}}
\author{Sayan Sadhukhan\thanks{Department of Computer Science} \and Warren Hare\thanks{Department of Mathematics} \and Yves Lucet\thanks{Department of Computer Science, corresponding author: yves.lucet@ubc.ca}\\I. K. Barber Faculty of Science, University of British Columbia, \\Okanagan Campus, Kelowna BC, Canada}

\begin{document}
	\maketitle
	\begin{abstract}
		The vertical alignment optimization problem in road design seeks the optimal vertical alignment of a road at minimal cost, taking into account earthwork while meeting all safety and design requirements.  In recent years, modelling techniques have been advanced to incorporate: side slopes, multiple material types, multiple hauling types, and road networks.  However, the advancements were created disjointly with implementations that only made a single advancement to the basic model.  Herein, we present a mixed-integer linear programming optimization model that unifies all previous advancements. 	The model further improves on previous work by maintaining convexity even in the multi-material setting. We compare our new model to previous models, validate it numerically, and demonstrate its capability in approximating material volumes. Our new model performs particularly well for determining the optimal vertical alignment for large road networks.
	\end{abstract}


\section{Introduction}\label{s:intro}

Transportation networks, encompassing roads and railway systems, are crucial in stimulating economic advancement, streamlining commerce, and guaranteeing the delivery of critical services \cite{RODRIGUE-06}. The task of crafting  roads is intricate and demands considerable time. Alignment Optimization has garnered attention, employing cutting-edge algorithmic solutions for thorough design and assessment of practical layouts, while adhering to an array of considerations, spanning safety, socio-economic factors, and environmental stipulations. The literature review paper in the field of alignment optimization~\cite{SONG-23} summarizes the accomplishments of the past 25 years.

The foundational aspect of a road design problem involves choosing a path that links the desired locations. This is inherently a multi-tiered mathematical optimization problem within the field of operations research and can be segmented into three distinct yet interrelated phases: horizontal alignment, vertical alignment, and earthwork procedures.

The process of road design initiates with the determination of a horizontal alignment, focusing on optimizing the trajectory of the road from an aerial perspective. Once the horizontal alignment is established, the focus shifts to the vertical alignment. This stage involves adapting a road profile to match the terrain, with an emphasis on minimizing construction costs while adhering to various constraints and road standards. This phase culminates in the derivation of a road profile as a solution. After fixing the vertical alignment, the attention turns to optimizing earthwork. This involves modeling construction site logistics to facilitate the terrain's transformation into a uniform surface, all the while keeping costs to a minimum.

\subsection{Horizontal alignment optimization}

Optimization of horizontal Alignment can be segmented into two distinct stages. Initially, the focus is on identifying an appropriate corridor for the intended road construction. Hirpa et al. \cite{HIRPA-16} present a biobjective optimization method for three-dimensional road alignment design while Pushak et al.~\cite{PUSHAK-16} propose a multiple path selection algorithm using adapted forms of Dijkstra's algorithm. An improved strategy building on that work allows to even consider the number of points of intersection as variables~\cite{ZHANG-23}.
A novel strategy aiming at optimal corridor selection with an emphasis on minimizing air pollution is discussed in \cite{GARCIA-CHAN-21}.

After determining the primary corridor, the process advances to the optimization of horizontal alignment within the designated corridor. This phase seeks to optimize the cost of construction and adhere to other constraints, encompassing diverse socioeconomic and political factors. Horizontal alignment optimization has been carried out using techniques like bi-level optimization~\cite{MONDAL-15},  network optimization~\cite{JHA-06}, and sequential quadratic programming~\cite{CASAL-17}.  Methodologies for optimizing horizontal alignment generally have highly complex models, approached by heuristics, or require repeatedly solving many vertical alignment and earthwork optimization problems.  This paper advances vertical alignment and earthwork optimization, which we discuss next.




\subsection{Vertical alignment and earthwork optimization}
Once we fix the horizontal alignment of the road,  the next step is to compute the vertical alignment. This phase centers on tailoring the road profile to match the ground's contour while minimizing costs and respecting various constraints and road specifications. Some of these constraints for vertical alignment include but are not limited to road grades and passenger safety. 

Once the vertical alignment of a road has been fixed, the earthwork optimization model is responsible for modeling the construction site logistics, ensuring that the earth is rearranged into a smooth surface based on the output of the vertical alignment model, while minimizing costs. At this stage,  optimization techniques can be used to optimize earthwork allocations, equipment fleet scheduling, and equipment routing. Earthwork optimization problems related to earthmoving operations can thus be classified into two further sub-categories, which are Equipment Fleet Planning and Earth Allocation Planning~\cite{GWAK-18}. 


Equipment Fleet Planning models revolve around picking the optimal ensemble of machinery, factoring in the machinery's kind and its operational efficiency across various excavation scenarios. Essentially, they sculpt a near-optimal team composition and an economical operational strategy.  Conversely, Earth Allocation Planning models cater to the blueprinting phase of a road design venture. These blueprints are instrumental in determining optimal earthwork distribution between excavation and embankment zones, also accounting for external earth sources, known as borrow pits, and designated zones for surplus material disposal, termed waste pits. A prevalent trend in academia is the unification of vertical alignment optimization and earthwork optimization into a single stage.

Stark and Nicholls \cite{STARK-72} were among the first people to show that linear programming models were superior in efficiency compared to traditional mass diagrams when orchestrating earthmoving logistics, leading to a more economical solution. Mayer and Stark \cite{MAYER-81} proposed a linear programming model for minimizing costs that included three categories of costs for excavation and loading, haul, and embankment, as well as included constraints for borrow pits, disposal sites, and shrinkage and swelling coefficients. In a different approach, Easa \cite{EASA-88} used a linear programming model for vertical alignment, which minimizes the total earthwork cost for finding the vertical alignment. Building upon previous works, Moreb \cite{MOREB-96} amalgamated both stages into a single model, ensuring global optimality, albeit at the expense of generating non-smooth solutions. To enhance this, Moreb and Aljohani \cite{MOREB-04} devised a road profile model utilizing quadratic functions, that removed the sharp connected ends of the linear segments.  

Koch and Lucet \cite{KOCH-10} refined the gap and slope constraints of the model, delivering increased precision in the solutions. They established that the model`s linearity persisted up to a quadratic spline. Hare et al.~\cite{HARE-11} introduced a mixed-integer linear programming approach tailored to the problem, which considered obstacles like rivers and mountain ranges. Later, Hare et al.~\cite{HARE-15} implemented a quasi-network flow (QNF) framework, building upon previous findings and introducing constraints for side slopes. Hossain \cite{HOSSAIN-13} extended this model to include multiple materials for earthwork allocation. Beiranvand et al.~\cite{BEIRANVAND-17} further expanded the model to integrate multiple haul routes, employing a multi-haul quasi-network flow (MH-QNF) approach for vertical alignment optimization.  Ayman et al.~\cite{AYMAN-23} further devised a mixed-integer linear programming solution targeting multiple interconnected roadways within a network using a multi-road quasi-network flow (MR-QNF) model. This approach, however, overlooked aspects like multiple materials, diverse hauling routes, and side slopes. Momo et al.~\cite{MOMO-23} proposed a quasi-network flow model using quadratic constraints. This gave rise to a convex model, termed the QCQP-QNF model, which catered to side slopes, but did not encompass considerations for multiple roads, diverse materials, or varied hauling paths. 

Our goal is to aggregate the above contributions in a unified model, and to propose a convex model that can be solved efficiently to global optimality.

The aforementioned techniques are inherently deterministic. However, there are nondeterministic strategies for determining vertical alignment documented in the academic literature. Lee and Cheng \cite{LEE-01} introduced a heuristic strategy that utilized multiple layers to transform a mixed-integer nonlinear program into its linear counterpart. Both Fwa et al. \cite{FWA-02} and Goektepe et al. \cite{GOKTEPE-09} employed genetic algorithms to address the vertical alignment challenge. For an in-depth exploration of alignment optimization solutions via genetic algorithms, one can refer to \cite{JHA-06}. In a separate vein, Goh et al. \cite{GOH-88} harnessed dynamic programming to model the optimal vertical alignment of a road while considering several costs such as cut, fill, pavement, and land acquisition costs.

\subsection{Three dimensional alignment optimization}\label{ss:3D}
In an ideal scenario, we aim to concurrently optimize horizontal alignment, vertical alignment, and earthwork allocation, a process termed as three-dimensional alignment optimization. Yet, practically speaking, this presents a significant challenge due to its non-convex nature. Most academic studies centered on 3D alignment optimization tend to employ nondeterministic algorithms with no guarantee of optimality  \cite{CHEW-89,TAT-03,LI-13,GAO-22}.

\subsection{Our research overview}

We focus on the optimization of vertical alignment and earthwork for large road networks within an acceptable computational time. While computing a vertical alignment is a necessary step to computing a horizontal alignment, it is also done routinely in practice when existing roads are upgraded by only optimizing the vertical alignment without changing the horizontal alignment. 

We adopt a convex model that scales in polynomial time.  This improves the previous MILP approach, which scales in exponential time.
Our work builds on the quasi network flow model \cite{HARE-11} for solving the vertical alignment problem and accounts for side slopes. We incorporate ideas from \cite{BEIRANVAND-17} that factor in multiple hauling paths and multiple materials, and \cite{AYMAN-23} which focuses on the modeling of multiple roads in a road network. The culmination of these enhancements results in the unified vertical alignment model named UVA detailed in Section~\ref{s:UVA}. 

In Section~\ref{s:CUVA}, we adapt the QCQP-QNF model from \cite{MOMO-23} to fit our multi-road multi-haul (UVA) framework, by translating our MILP volume constraints into quadratic constraints. We extend this approach to also factor in the volume of multiple materials. This leads to a new model called the convex unified vertical alignment (CUVA) model.

Section~\ref{s:val} is dedicated to the validation of our CUVA, where we utilize data from \cite{MOMO-23,BEIRANVAND-17,AYMAN-23}. We show that our CUVA model does a better job at approximating volumes. Lastly, Section~\ref{s:conc} provides avenues for subsequent research.

\section{QNF model variables and parameters}\label{s:UVA}

The multi-road multi-haul quasi-network flow (QNF) model aggregates the quasi network flow model \cite{HARE-11} for solving the vertical alignment problem while taking side slopes into account with the multi-haul multi-material model from~\cite{BEIRANVAND-17}, and the multi-road model from~\cite{AYMAN-23}. We now outline its variables (colored in \var{red}) and its parameters.

Consider a road network with $r$ roads indexed by $\mathcal{R} = \{1,2,3,...,r\}$. Each $i^{th}$ road vertical profile is a quadratic spline split into $n_{seg_{i}}$ segments where each segment is indexed by $\mathcal{G}_{i} = \{1,2,3,...,n_{seg_{i}}\}$. Each $g^{th}$ spline segment of the $i^{th}$ road is further comprised of $n_{i,g}$ sections indexed by $\mathcal{S}_{i,g} = \{1,2,3,...,n_{i,g}\}$. The total number of sections in a road $i$ is $n_{i}=\sum_{g \in \mathcal{G}_{i}} n_{i,g}$. We denote $\varphi : (\mathcal{R}\times\mathcal{G}_{i}\times\mathcal{S}_{i,g})\rightarrow \mathcal{S}_{i}$ the function that maps the section indexes of a particular spline segment $g_{i}$ to the actual section index set $\mathcal{S}_{i}$ for the $i^{th}$ road. For example, if $\varphi(i,g,j)=s $ then $x_{s} = x_{\varphi(i,g,j)}$ for all $s \in \mathcal{S}_{i},g \in \mathcal{G}_{i},j \in \mathcal{S}_{i,g}, i \in \mathcal{R}$. 

\begin{figure}
	\centering
	\includegraphics[width=0.8\textwidth]{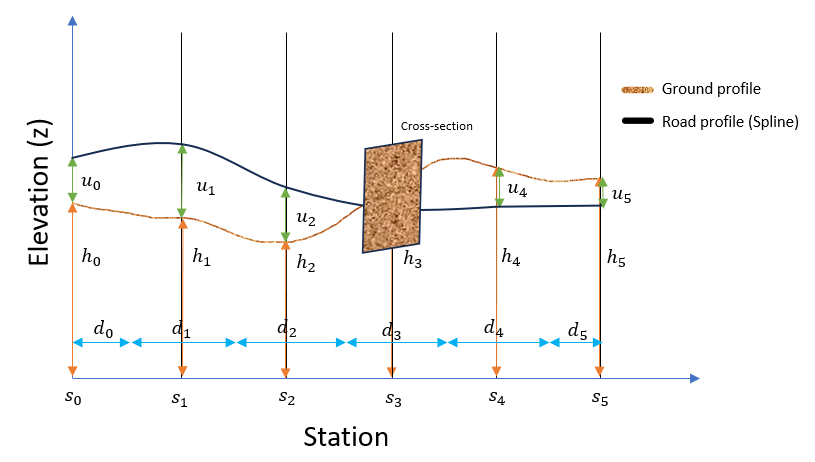}
	\caption{Side view of the ground profile and road profile.}
	\label{fig:vert-align-road}
\end{figure}

Every pair of sequential segments share a mutual station, see Figure~\ref{fig:vert-align-road}. In this context, $d_{i}$ represents the effective length for a section, $h_{i}$ denotes the elevation of the ground, while $u_{i}$ is the difference in elevation between the ground profile and the road profile. The terminal station of one segment coincides with the starting station of its subsequent segment. At every station, a snapshot of the terrain's cross-section is taken. This assists in determining the cross-sectional surface area for every type of material identified. By scaling this area with the effective length of the segment, we can approximate the volume of material. While a section is associated with material volume, a station focuses on the cross-sectional surfaces at a designated point. The QNF model handles non-uniform station locations, however engineers typically use a constant section length. 

The altitude of a point at distance $s$ from the beginning of Segment $g$ on Road $i$ is given by
\begin{equation}
	\label{eqn:SplineEqQuadratic}
	\var{P_{i,g}}(s) = \var{a_{i,g,1}} + \var{a_{i,g,2}} (s - s_{i,g}) + \var{a_{i,g,3}} (s - s_{i,g})^2.
\end{equation}
%
We denote $\var{P'_{i}}(s)$ its derivative.

Within a road network, intersections are points where multiple roads meet.  Consider $\mathcal{I} = \{1,2,3,...,n_{e}\}$ as the index set representing all intersections within the road network, and $\mathcal{R}_{e} = \{1,2,3,...,n_{r,e}\}$ as the index set denoting all roads linked to intersection $e$. The function $\xi :\mathcal{I} \rightarrow \mathcal{R}_{e}$ yields the set of roads joining intersection $e$, while the function $\gamma: (\mathcal{R}_{e}\times\mathcal{I})\rightarrow\mathcal{S}_{i,e}$ provides the set of sections from road $i$ linked to intersection $e$. Furthermore, the functions $\mu :(\mathcal{I}\times\mathcal{R}_{e})\rightarrow\mathcal{R}$ and $\tau :(\mathcal{I}\times\mathcal{R}_{e})\rightarrow\mathcal{S}_{i}$ indicate the road and section indices associated with a particular intersection. In a similar vein, the functions $\eta :(\mathcal{I}\times\mathcal{R}_{e})\rightarrow\mathcal{G}_{i}$ and $\zeta :(\mathcal{I}\times\mathcal{R}_{e})\rightarrow\mathcal{S}_{i,g}$ map the segment and section connected to an intersection to their respective index sets, $\mathcal{G}_{i}$ and $\mathcal{S}_{i,g}$. Lastly, the function $\lambda : (\mathcal{R}\times\mathcal{S}_{i})\rightarrow\mathcal{I}$ identifies the intersection linked to the $i^{th}$ road's $j^{th}$ section, where $i \in \mathcal{R}$ and $j \in \mathcal{S}_{i}$.

In a typical real-world road construction project various materials are used, each with its own associated costs for excavation, embankment, and hauling. Multiple material types are accounted for in our optimization model to generate more accurate results. Let the index set of the material types be denoted by $\mathcal{K}=\{1,2,\dots,m\}$.

In the vertical alignment problem, we strive to obtain the optimal offset between the ground profile and the road profile. We also need to calculate the required cut and fill volumes for each offset. The offset is denoted by $\var{u_{i,j}}$, for each $i \in \mathcal{R}$ and for each $j \in \mathcal{S}$, and the cut and fill volumes of a section $j$ are denoted by $\var{V_{i,j,m}^{+}}$ and $\var{V_{i,j,m}^{-}}$ respectively, for each road $i \in \mathcal{R}$, section $j \in \mathcal{S}$, and material $m \in \mathcal{K}$. It is important to highlight that the boundaries for $\var{u_{i,j}}$ are set by $\underline{u}_{i,j}$ for the lower limit and $\overline{u}_{i,j}$ for the upper limit. In this context, a negative $\var{u_{i,j}}$ indicates an excavation, a positive $\var{u_{i,j}}$ indicates a fill, and when $\var{u_{i,j}}$ equals zero, it implies that a section has not been subjected to either excavation or fill. 

In the absence of side-slopes, the volume of section $j$ for road $i$ containing $m$ materials can be described as a linear function of $\var{u_{i,j}}$. However, when side slopes are incorporated, the road's cross-section at a particular section assumes a trapezoidal shape. As a result, the volume evolves into a quadratic function in terms of $\var{u_{i,j}}$, which disrupts the model's linear nature. Yet, as highlighted in~\cite{HARE-15}, the model's linearity is preserved by approximating the excavation volume of section $j$ for road $i$ with $m$ materials using $k^{+}_{i,j,m}$ rectangular sections, and the fill volume using $k^{-}_{i,j,m}$ rectangular sections. These rectangular sections are termed ``slabs''. 
For an excavated section $j$ of road $i$ with material $m$, the slab areas and heights are denoted by $A_{i,j,m}^{k+}$ and $\overline{u}_{i,j,m}^{k+}$ respectively, where $k$ belongs to the set $\{1,2,\dots,k_{i,j,m}^{+}\}$. For fill, these are represented as $A_{i,j,m}^{k-}$ and $\underline{u}_{i,j,m}^{k-}$ with $k$ in the set $\{1,2,\dots,k_{i,j,m}^{-}\}$. 
Let the parameter $\mathcal{P}_{i,j,m_1}^{+}$ represent the cut and $\mathcal{P}_{i,j,m_1}^{-}$ represent the fill volume for each road $i \in \mathcal{R}$, each section $j \in \mathcal{S}$, and for material $m_1 \in \mathcal{K}$.

%
%

Two sets of indices $\mathcal{B}=\{1,2,3,...,n_{i,\beta}\}$ and $\mathcal{W}=\{1,2,3,...,n_{i,\omega}\}$ are used to index borrow pits and waste pits respectively, where $n_{i,\beta}$ is the number of borrow pit sections and $n_{i,\omega}$ is the number of waste pit sections for the $i^{th}$ road. We map the borrow pit index to the corresponding section index to which it is attached using $\nu: \mathcal{B}_{i}\rightarrow\mathcal{S}_{i}$. Similarly, for mapping the waste pit index to its corresponding section index we use the function $\delta: \mathcal{W}_{i}\rightarrow\mathcal{S}_{i}$. The capacity of a borrow pit $j$ for the $i^{th}$ road is denoted by $\mathcal{C}_{i,j,m}^{b}$, and the capacity of a waste pit $k$ for the $i^{th}$ road is denoted by $\mathcal{C}_{i,k,m}^{w}$ for each material $m \in \mathcal{K}$. The distance between the $j^{th}$ borrow/waste pit and the section it is connected to is called dead haul distance and is represented by $\widetilde{d_{j}}$ and $\widetilde{d_{k}}$ for borrow and waste pit respectively, where $j \in \mathcal{B}_{i}$ and $k \in \mathcal{W}_{i}$ for road $i$. The set $\mathcal{N}_{i} = \{\mathcal{S}_{i}\cup\mathcal{B}_{i}\cup\mathcal{W}_{i}\}$ is used to collectively represent all the indexes for sections, borrow pits, and waste pits respectively for road $i$.

\section{Unified Vertical Alignment model (UVA)} 

The UVA model extends the basic QNF model to include intersections, and hauling paths. Initially, we delve into the single-haul scenario, where different materials are transported via a single hauling route, before transitioning to a more comprehensive multi-haul strategy. Analogous to the original QNF framework \cite{HARE-15}, we represent sections, crossroads, borrow pits, and waste pits as nodes. Additionally, we depict viable routes as arcs to facilitate the transfer of materials between nodes. 


\subsection{Model Description}


In the initial QNF framework, the researchers employed a single haul approach for material movement. This refers to using only one kind of earthmoving machinery. In projects related to road construction, various types of such machinery might be utilized for distinct purposes. For instance, bulldozers (resp. front loader, truck) might be economically efficient for short (resp. medium, long) spans. To make the model more representative of real-world scenarios, the original QNF was enhanced to include multiple hauling routes \cite{BEIRANVAND-17}.


The objective of the model is to ascertain the optimal hauling path or the most cost-effective equipment based on the hauling distance. Within this framework, materials can be transported from the present section $j$ either to the preceding section $(j - 1)$ or the subsequent section $(j + 1)$, for road $i$. The transfer of material can be executed via any of the hauling paths denoted by $\mathcal{H} = \{1, 2, ..., n_{h}\}$. For our case, we shall apply $n_{h} = 3$ for our test problems. 

Upon arriving at a section, the material has three possible actions: it can be offloaded to the section, extracted from the section introducing new material into the flow, or simply passed onto the subsequent section without any change. We employ virtual nodes to facilitate material transfer in both leftward and rightward directions. With $n_{h}$ distinct hauling paths, the model establishes $n_{h}$ groups of transition nodes.  

For example, if $n_{h}=3$, then we have $3$ distinct hauling paths, which we shall refer to as short-haul, middle-haul, and long-haul.
In this case, we define variables $\var{f_{i,j,m}^{t^+,h1}}$ and $\var{f_{i,j-1,m}^{t^+,h1}}$ (for all $j \in \mathcal{S}$) as the flow of material from the left (or right) short-haul transit node to the right (or left) short-haul transit node. Similarly, variables $\var{f_{i,j,m}^{t^+,h2}}$ and $\var{f_{i,j-1,m}^{t^+,h2}}$ (for all $j \in \mathcal{S}$) are defined for the middle-haul path, and $\var{f_{i,j,m}^{t^+,h3}}$ and $\var{f_{i,j-1,m}^{t^+,h3}}$ (for all $j \in \mathcal{S}$) denote the variables for the long-haul path. 

To transfer material from section $j$ to section $j+1$, it is necessary to extract material from section $j$. This extracted material is then loaded onto the appropriate equipment, considering the hauling distance (essentially, the material is taken out of section $j$). This loading task is modeled using loading flow variables $\var{f_{i,j,m}^{l^+,h}}$ and $\var{f_{i,j-1,m}^{l^-,h}}$ (for all $j \in \mathcal{S}$, $h \in \mathcal{H}$). Likewise, for the target section $j$, the process of depositing material into the section is represented using the unloading flow variables $\var{f_{i,j,m}^{u^+,h}}$ and $\var{f_{i,j-1,m}^{u^-,h}}$ (for all $j \in \mathcal{S}$, $h \in \mathcal{H}$). Once the materials have been excavated, they can be transported via any available hauling methods. 

In a comparable manner, a section can be replenished with materials sourced from either the left or right transition nodes. To model this scenario, we use the variables $\var{f_{i,j,m}^{l^+,h}}$ and $\var{f_{i,j-1,m}^{l^-,h}}$ (for all $j \in \mathcal{S}$) as the load flows of materials from the left and the right transit nodes for the hauling path $h$. Similarly, for the unload flows of materials to the left and right nodes we define, $\var{f_{i,j,m}^{u^+,h}}$, $\var{f_{i,j-1,m}^{u^-,h}}$ (for all $j \in \mathcal{S}$).


The variables $\var{f^{b,h}_{i,\vartheta(k)-1,m}}$ and $\var{f^{b,h}_{i,\vartheta(k),m}}$ (for all $k \in \mathcal{B}$) represent the borrow flows from the borrow pit j to both left and right directions through hauling path $h$. The variables $\var{f^{w,h}_{i,\delta(l)-1,m}}$ and $\var{f^{w,h}_{i,\delta(l),m}}$ (for all $l \in \mathcal{W}$) represent the waste flows to the waste pit $k$ from both left and right directions through hauling path $h$. 

Contrary to the QNF model, the UVA model is capable of managing multiple roads and their corresponding intersections. 
%
%

Within this model, every road in the network is segmented into distinct sections. Each of these sections boasts two virtual transit nodes: a positive node for facilitating material movement from left to right, and a negative one for directing material flow from right to left. In other words, for each road, we have a start section and an end section. The positive nodes facilitate the movement of materials from the start section to the end section, while the negative nodes facilitate the movement of materials in the reverse direction, that is from the end section to the start section. The positive transit node possesses an outgoing flow path, represented as $\var{f_{i,j,m}^{t+,h}}$ (where $i \in \mathcal{R}, j \in \mathcal{S}_{i}, m \in \mathcal{K}, h \in \mathcal{H}$), except for the terminal section's node. Additionally, it has an incoming flow path, $\var{f_{i,j-1,m}^{t+,h}}$, excluding the inaugural section's node. In a parallel fashion, the negative transit node features an outgoing flow path, symbolized by $\var{f_{i,j,m}^{t-,h}}$ (with $i \in \mathcal{R}, j \in \mathcal{S}_{i}, m \in \mathcal{K}, h \in \mathcal{H}$), barring the first section's node, and an inbound flow path, $\var{f_{i,j-1,m}^{t-,h}}$, except for the node of the final section. Furthermore, both types of nodes (positive and negative) are equipped with an additional pair of flow paths designated for material loading and unloading tasks. Here, $\var{f_{i,j,m}^{l+,h}}$ and $\var{f_{i,j,m}^{l-,h}}$ represent the dual outbound loading flow paths, while $\var{f_{i,j,m}^{u+,h}}$ and $\var{f_{i,j,m}^{u-,h}}$ denote the pair of inbound unloading flow paths for the positive and negative transit nodes respectively.  If a specific section is linked to an access road, enabling material transport to either a borrow pit or a waste pit, then its transit nodes gain two more flow paths. This includes the outbound paths $\var{f_{i,j,m}^{w+,h}}$ and $\var{f_{i,j,m}^{w-,h}}$ for directing excess materials to a waste pit and the inbound paths $\var{f_{i,j,m}^{b+,h}}$ and $\var{f_{i,j,m}^{b-,h}}$ for sourcing materials from a borrow pit.

An intersection represents a shared section between multiple roads. 
Each road can intersect with a maximum of two other roads within the scope of the model. Specifically, a road can only intersect at its beginning or end. For instance, if two roads cross at a shared section, the model interprets the road configuration as four distinct roads converging at that shared section or intersection.
Every intersection in this model is equipped with a single virtual transitional node, responsible for directing materials between roads. This virtual node encompasses two flow pathways for both loading and offloading materials onto the intersection. The symbols $\var{f_{e,m}^{l,h}}$ and $\var{f_{e,m}^{u,h}}$ stand for the outbound loading route and inbound unloading route, respectively, for the intersection's virtual transitional node, where $e$ belongs to set $\mathcal{I}$, and $m$ belongs to $\mathcal{K}$, and $h$ belongs to $\mathcal{H}$. Additionally, the virtual node at an intersection is endowed with two supplementary flow paths (one for intake and one for outflow) for each road that it intersects. 
The notations $\var{f_{i,j,m}^{in,e,h}}$ and $\var{f_{i,j,m}^{out,e,h}}$ represent the respective inbound and outbound flow paths for every intersecting road.

\subsection{Cost Parameters}

The UVA model takes into account the costs of excavation, embankment, hauling, and loading materials for the construction of a road network. Let the symbols $p_1, p_2, \dots p_m$ represent the excavation costs per unit volume for various materials. Similarly, $q_1, q_2, \dots q_m$ stand for the embankment costs per unit volume for these different materials. 

For every road and hauling path, there exist two potential directions for material movement. For hauling path $h$, the hauling cost of materials from section $j$ to section $j+1$ is defined as $c_{i,j,m}^{t^+,h}$ = $c_{h,m}d_{i,j^+}$, and the hauling cost of materials from section $j$ to section $j-1$ is defined as $c_{i,j-1,m}^{t^-,h}$ = $c_{h,m}d_{i,j^-}$. 

An additional expense in extracting materials from a section or borrow pit pertains to the loading cost. Practically, distinct construction equipment incurs varying loading expenses. As such, we characterize the loading costs for different equipment types using short, middle, and long hauling paths as $y_{h_{1},m}$, $y_{h_{2},m}$, and $y_{h_{3},m}$ for every material classified under $m \in \mathcal{K}$. As recommended in \cite{BEIRANVAND-17}, this model omits unloading costs, deeming them negligible in comparison to other expenses.

\subsection{Objective Function}

The objective of the UVA model is to minimize the total construction cost which includes excavation, embankment, hauling, and loading costs of materials for constructing an entire road network. Hence, the objective function to minimize is
\begin{multline} \label{eqn:objectiveFunc}
		  \smashoperator{\sum\limits_{
				\begin{subarray}{c}
					i \in \mathcal{R} \\
					j \in \mathcal{S}_i \cup \mathcal{B}_i \\
					m \in \mathcal{K} \\
					h \in \mathcal{H}
				\end{subarray}
		}}
		\left(p_{m,h}+y_{m,h}\right) \var{\mathcal{V}_{i,j,m,h}^{+}} + \smashoperator{\sum\limits_{
				\begin{subarray}{c}
					i \in \mathcal{R} \\
					j \in \mathcal{S}_i  \cup \mathcal{W}_i \\
					m \in \mathcal{K} \\
					h \in \mathcal{H}
				\end{subarray}
		}} q_c\var{\mathcal{V}_{i,j,m,h}^{-}} 
		 + \smashoperator{\sum\limits_{
				\begin{subarray}{c}
					i \in \mathcal{R} \\
					j \in \mathcal{S}_i \\
					m \in \mathcal{K} \\
					h \in \mathcal{H}
				\end{subarray}
		}}  \left(c_{i,j,m}^{t+,h} \var{f_{i,j,m}^{t+,h}} + c_{i,j,m}^{t-,h} \var{f_{i,j,m}^{t-,h}} \right) \\
		 + \smashoperator{\sum\limits_{
				\begin{subarray}{c}
					e \in \mathcal{I} \\
					m \in \mathcal{K} \\
					h \in \mathcal{H}
				\end{subarray}
		}} \left(p_{m,h}+y_{m,h}\right) \var{\mathcal{U}_{e,m,h}^{+}} + \smashoperator{\sum\limits_{
				\begin{subarray}{c}
					e \in \mathcal{I} \\
					m \in \mathcal{K} \\
					h \in \mathcal{H}
				\end{subarray}
		}} q_{m,h} \ \var{\mathcal{U}_{e,m,h}^{-}} 
		 + \smashoperator{\sum\limits_{
				\begin{subarray}{c}
					e \in \mathcal{R} \\
					i \in \xi(e) \\
					j \in \gamma(i,e) \\
					m \in \mathcal{K} \\
					h \in \mathcal{H}
				\end{subarray}
		}} c_{i,j,m}^{e,h}\left( \var{f_{i, j, m}^{in, e, h}} + \var{f_{i, j, m}^{out, e, h}} \right) \\
		 + \smashoperator{\sum\limits_{
				\begin{subarray}{c}
					i \in \mathcal{R} \\
					k \in \mathcal{B}_{i}\\
					m \in \mathcal{K} \\
					h \in \mathcal{H}
				\end{subarray}
		}} \left(p_{m,h} + y_{m,h} + c_{m,h} \tilde{d}_{k}\right) \left(\var{f_{i, \vartheta(k),m}^{b+,k,h}} + \var{f_{i, \vartheta(k),m}^{b-,k,h}}\right) 
		 + \smashoperator{\sum\limits_{
				\begin{subarray}{c}
					i \in \mathcal{R} \\
					l \in \mathcal{W}_i \\
					m \in \mathcal{K} \\
					h \in \mathcal{H}
				\end{subarray}
		}} \left(q_{m,h} + c_{m,h} \tilde{d}_{l}\right) \left(\var{f_{i, \varphi(l), m}^{w+,l,h}} + \var{f_{i, \varphi(l),m}^{w-,l,h}}\right).
\end{multline}

\subsection{Constraints}
\begin{description}[style=unboxed,leftmargin=0cm]
\item[Continuity constraints]
The elevation and gradient of a segment's initial section must align with the elevation and gradient of the preceding segment's final section, i.e., 
\begin{equation}
	\label{eqn:ContinuityConsRoad1}
	\begin{aligned}
		\var{P_{i,g-1}}\left(x_{\varphi \left(i,g-1,n_{i,{g-1}}\right)}\right) = \var{P_{i,g}}\left(x_{\varphi \left(i,g,1\right)}\right) && \forall i \in \mathcal{R}, \forall g \in \mathcal{G}_i \setminus \{1\}, 
	\end{aligned}
\end{equation}
\begin{equation}
	\label{eqn:ContinuityConsRoad2}
	\begin{aligned}
		\var{P^\prime_{i,g-1}}\left(x_{\varphi \left(i,g-1,n_{i,{g-1}}\right)}\right) = \var{P^\prime_{i,g}}\left(x_{\varphi \left(i,g,1\right)}\right) && \forall i \in \mathcal{R}, \forall g \in \mathcal{G}_i \setminus \{1\}. 
	\end{aligned}
\end{equation}

Likewise, the elevation of all sections converging at a single intersection should be consistent
\begin{multline}
	\label{eqn:ContinuityConsIntersection}
	\var{P_{\mu(e,1),\eta(e,1)}}\left(x_{\varphi \left(\mu(e,1),\eta(e,1),\zeta(e,1)\right)}\right) = 
	\var{P_{\mu(e,n_{re}),\eta(e,n_{re})}}\left(x_{\varphi \left(\mu(e,n_{re}),\eta(e,n_{re}),\zeta(e,n_{re})\right)}\right) \\
	\forall e \in \mathcal{I}, \forall i \in \mathcal{R}.
\end{multline}

	
\item[Gap constraints]
The gap constraints ascertain that the variance between the road's profile and the terrain's profile for every section within a segment aligns with that section's offset. For all $i \in \mathcal{R}, j \in \mathcal{S}_i, g \in \mathcal{G}_i, k \in \mathcal{S}_{i,g},$ the gap constraint is defined as

\begin{equation}
	\label{eqn:GapCons1}
	\begin{aligned}
		\var{P_{i,g}}\left(x_{\varphi\left(i,g,k\right)}\right) - Z_{i,j} = \var{u_{i,j}}.
	\end{aligned}
\end{equation}

Here, $Z_{i,j}$ is the altitude of the terrain for $i^{th}$ road's $j^{th}$ section. For an intersection, the gap constraint is defined as
\begin{equation}
	\label{eqn:GapCons2}
	\begin{aligned}
		\var{u_{e}} = \var{u_{\xi(e)}}, \forall e \in \mathcal{I}.
	\end{aligned}
\end{equation}
This means that for an intersection, its offset must be consistent with the offset of the adjacent section from various roads to which it is linked.

\item[Grade constraints]
The gradient constraints are instituted to uphold road safety. They confine the road gradients within a designated range, bounded by a minimum value, $G_L$, and a maximum value, $G_U$. The road profile is constructed adhering to this spectrum. For all $i \in \mathcal{R}, g \in \mathcal{G}_i$, the grade constraints are defined as 

\begin{equation}
	\label{eqn:GradeConsRoad1}
	\begin{aligned}
		G_L \leq \var{P^\prime_{i,g}}\left(x_{\varphi \left(i,g,1\right)}\right) \leq G_U,
	\end{aligned}
\end{equation}
\begin{equation}
	\label{eqn:GradeConsRoad2}
	\begin{aligned}
		G_L \leq \var{P^\prime_{i,g,n_g}}\left(x_{\varphi \left(i,g,n_g\right)}\right) \leq G_U.
	\end{aligned}
\end{equation}

\item[Flow constraints]
Flow constraints guarantee that the cumulative quantity of materials entering a specific virtual transit node of a section is equivalent to the total materials exiting that node. The flow constraints for virtual transit nodes when the section is not linked to an intersection, are for all $i \in \mathcal{R}, j \in \mathcal{S}_i, m \in \mathcal{K}, h \in \mathcal{H},$
\begin{align}
		\var{f^{t+,h}_{i,j,m}} + \var{f^{l+,h}_{i,j,m}} + \sum\limits_{l \in \mathcal{W}_i}\var{f^{w+,l,h}_{i,\delta(l),m}} &= \var{f^{t+,h}_{i,j-1,m}} + \var{f^{u+,h}_{i,j,m}}  + \sum\limits_{k \in \mathcal{B}_i}\var{f^{b+,k,h}_{i,\nu(k),m}},
		\label{eqn:FlowConsRoad1} \\
		\var{f^{t-,h}_{i,j,m}} + \var{f^{l-,h}_{i,j,m}} + \sum\limits_{l \in \mathcal{W}_i}\var{f^{w-,l}_{i,\delta(l)}} &= \var{f^{t-,h}_{i,j-1,m}} + \var{f^{u-,h}_{i,j,m}} + \sum\limits_{k \in \mathcal{B}_i}\var{f^{b-,k,h}_{i,\nu(k),m}}.\label{eqn:FlowConsRoad2}
\end{align}

In scenarios where a specific section is linked to an intersection and represents the initial section of the road, the flow constraints for the transit nodes of this section are delineated as
\begin{align}
		\var{f^{t+,h}_{i,j,m}}  + \sum\limits_{l \in \mathcal{W}_i}\var{f^{w+,l,h}_{i,\delta(l)}} &= \var{f^{out,\lambda(i,j),h}_{i,j,m}} + \sum\limits_{k \in \mathcal{B}_i}\var{f^{b+,k,h}_{i,\nu(k),m}}, 	\label{eqn:FlowConsRoad3} \\
		\var{f^{in,\lambda(i,j),h}_{i,j,m}} + \sum\limits_{l \in \mathcal{W}_i} \var{f^{w-,m,h}_{i,\delta(m),m}} &= \var{f^{t-,h}_{i,j-1,m}} + \sum\limits_{k \in \mathcal{B}_i}\var{f^{b-,k,h}_{i,\nu(k),m}}.
	\label{eqn:FlowConsRoad4}
\end{align}

Likewise, if a section, acting as the terminal segment of a road, is linked to an intersection, then the corresponding flow constraints are delineated as
\begin{align}
		\var{f^{in,\lambda(i,j),h}_{i,j,m}}  + \sum\limits_{l \in \mathcal{W}_i} \var{f^{w+,l,h}_{i,\delta(l),m}} &= \var{f^{t+,h}_{i,j-1,m}} + \sum\limits_{k \in \mathcal{B}_i}\var{f^{b+,k,h}_{i,\nu(k),m}},
	\label{eqn:FlowConsRoad5} \\
		\var{f^{t-,h}_{i,j,m}} + \sum\limits_{l \in \mathcal{W}_i}\var{f^{w-,l,h}_{i,\delta(l),m}} &= \var{f^{out,\lambda(i,j),h}_{i,j,m}} + \sum\limits_{k \in \mathcal{B}_i}\var{f^{b-,k,h}_{i,\nu(k),m}}.
		\label{eqn:FlowConsRoad6}
\end{align}

We further incorporate a constraint that maintains the material flow consistency at each intersection's transit node
\begin{equation}	\label{eqn:FlowConsInt1}
		\sum\limits_{i \in \xi(e)} \var{f^{in,e,h}_{i,\tau(e,i),m}} + \var{f^{u,h}_{e,m}} = \sum\limits_{i \in \xi(e)} \var{f^{out,e,h}_{i,\tau(e,i),m}} + \var{f^{l,h}_{e,m}},
\end{equation}
for all $e \in \mathcal{I}, m \in \mathcal{K}, h \in \mathcal{H}$.

Lastly, we establish constraints to guarantee that the loading and unloading of materials from sections linked to an intersection are not redundantly tallied, given that these materials are already accounted for within the intersection; for all $e \in \mathcal{I}, m \in \mathcal{K}, h \in \mathcal{H}$
\begin{align}
		\sum\limits_{i \in \xi(e)} \var{f^{u+,h}_{i,\tau(e,i),m}} + \sum\limits_{i \in \xi(e)}\var{f^{u-,h}_{i,\tau(e,i),m}} &= 0, 	\label{eqn:FlowConsInt2}\\	
		\sum\limits_{i \in \xi(e)} \var{f^{l+,h}_{i,\tau(e,i),m}} + \sum\limits_{i \in \xi(e)}\var{f^{l-,h}_{i,\tau(e,i),m}} &= 0. \label{eqn:FlowConsInt3}
\end{align}

\item[Balance constraints]
The cumulative unloading flows departing from a section and the total loading flows entering this section should correspond to the section's overall cut and fill volumes, respectively. These constraints are termed balance constraints and are defined for all $i \in \mathcal{R}, j \in \mathcal{S}_i, m \in \mathcal{K}$ as
\begin{align}
		\sum_{h \in \mathcal{H}} \var{f^{u+}_{i,j,m}} + 
		\sum_{h \in \mathcal{H}} \var{f^{u-}_{i,j,m}} &= 
		\sum_{h \in \mathcal{H}} \var{\mathcal{V}^{+}_{i,j,m}} , 	\label{eqn:BalConsRoad1} \\
		\sum_{h \in \mathcal{H}} \var{f^{l+}_{i,j,m}} + 
		\sum_{h \in \mathcal{H}} \var{f^{l-}_{i,j,m}} &= 
		\sum_{h \in \mathcal{H}} \var{\mathcal{V}^{-}_{i,j,m}}. 	\label{eqn:BalConsRoad2}
\end{align}

Similarly, for the intersection, the balance constraints are for all $e \in \mathcal{I}, m \in \mathcal{K}$
\begin{align}
		\sum_{h \in \mathcal{H}} \var{f^{u}_{e,m}} &= 
		\sum_{h \in \mathcal{H}} \var{\mathcal{U}^{+}_{e,m}}, \label{eqn:BalConsInt1} \\
		\sum_{h \in \mathcal{H}} \var{f^{l}_{e,m}} &= 
		\sum_{h \in \mathcal{H}} \var{\mathcal{U}^{-}_{e,m}}.	\label{eqn:BalConsInt2}
\end{align}

\item[Volume constraints]
Volume constraints ensure that, within a specific section, the cumulative volume either added to or extracted from the section matches the volume disparity between the road and terrain profiles. Given assumptions, $0<A_{i,j,m}^{1+}<A_{i,j,m}^{2+}<\dots<A_{i,j,m}^{\overline{k}_{i,j}^+}$ and $0<A_{i,j,m}^{1-}<A_{i,j,m}^{2-}<\dots<A_{i,j,m}^{\overline{k}_{i,j}^-}$, the volume constraints using the side-slope approximation are delineated for all $i \in \mathcal{R}, j \in \mathcal{S}_i, m \in \mathcal{K}$; we compute $\var{k_{i,j,m}^{+}} \in \{1,2,\dots,\overline{k}_{i,j,m}^{+}\}$, and $\var{k_{i,j,m}^{-}} \in \{1,2,\dots,\overline{k}_{i,j,m}^{-}\}$
as
\begin{align}
		\var{\mathcal{V}^+_{i,j,m}} &\geq \sum\limits_{k =1}^{\var{k_{i,j,m}^{+}} -1} A^{k}_{i,j,m} \overline{u}_{i,j,m}^{k} + A^{k^+}_{i,j,m} \left (\var{u_{i,j}} - \sum\limits_{k =1}^{\var{k_{i,j,m}^{+}} -1}\overline{u}_{i,j,m}^{k} \right ),
	\label{eqn:VolConsRoad1} \\
		\var{\mathcal{V}^-_{i,j,m}} &\geq - \sum\limits_{k =1}^{\var{k_{i,j,m}^{-}} -1} A^{k}_{i,j,m} \underline{u}_{i,j,m}^{k} - A^{k^-}_{i,j,m} \left (\var{u_{i,j}} - \sum\limits_{k =1}^{\var{k_{i,j,m}^{-}} -1}\underline{u}_{i,j,m}^{k} \right ). 	\label{eqn:VolConsRoad2}		
\end{align}

The cross-section of a cut section $j$ for road $i$ and material $m$ is approximated by $\overline{k}_{i,j}^+$ slabs having elevations $\overline{u}_{i,j,m}^{1},\overline{u}_{i,j,m}^{2},\dots,\overline{u}_{i,j,m}^{\overline{k}_{i,j}^+}$ and areas $A^{1+}_{i,j,m},A^{2+}_{i,j,m},\dots,A^{\overline{k}_{i,j}^+}_{i,j,m}$. Similarly, the cross-section of a fill section $j$ for road $i$ and material $m$ is approximated by $\overline{k}_{i,j,m}^{-}$ slabs having elevations $\underline{u}_{i,j,m}^{1},\underline{u}_{i,j,m}^{2},\dots,\underline{u}_{i,j,m}^{\overline{k}_{i,j,m}^{-}}$ and areas $A^{1-}_{i,j,m},A^{2-}_{i,j,m},\dots,A^{\overline{k}_{i,j,m}^{-}}_{i,j,m}$.

In a similar manner, the constraints for the intersection are defined for all
$e \in \mathcal{I}, m \in \mathcal{K}$ by computing $\var{k_{e,m}^{+}} \in \{1,2,\dots,\overline{k}_{e,m}^{+}\}$ and $\var{k_{e,m}^{-}} \in \{1,2,\dots,\overline{k}_{e,m}^{-}\}$
as
\begin{align}
		\var{\mathcal{U}^+_{e,m}} &\geq \sum\limits_{k =1}^{\var{k_{e,m}^{+}} -1} A^{k}_{e,m} \overline{z}_{e,m}^{k} + A^{k^+}_{e,m} \left (\var{z_{e}} - \sum\limits_{k =1}^{\var{k_{e,m}^{+}} -1}\overline{z}_{e,m}^{k} \right ),
		\label{eqn:VolConsInt1} \\
		\var{\mathcal{U}^-_{e,m}} &\geq - \sum\limits_{k =1}^{\var{k_{e,m}^{-}} -1} A^{k}_{e,m} \underline{z}_{e,m}^{k} - A^{k^-}_{e,m} \left (\var{z_{e}} - \sum\limits_{k =1}^{\var{k_{e,m}^{-}} -1}\underline{z}_{e,m}^{k} \right ). \label{eqn:VolConsInt2}
\end{align}

The cross-section of a cut section for intersection $e$ and material $m$ is approximated by $\overline{k}_{e,m}^{+}$ slabs having elevations $\overline{z}_{e,m}^{1},\overline{z}_{e,m}^{2},\dots,\overline{z}_{e,m}^{\overline{k}_{e,m}^{+}}$ and areas $A^{1+}_{e,m},A^{2+}_{e,m},\dots,A^{\overline{k}_{e,m}^{+}}_{e,m}$. Similarly, the cross-section of a fill section for intersection $e$ and material $m$ is approximated by $\overline{k}_{i,j,m}^{-}$ slabs having elevations $\underline{z}_{e,m}^{1},\underline{z}_{e,m}^{2},\dots,\underline{z}_{e,m}^{\overline{k}_{i,j,m}^{-}}$ and areas $A^{1-}_{e,m},A^{2-}_{e,m},\dots,A^{\overline{k}_{i,j,m}^{-}}_{e,m}$, with the assumptions $0<A_{e,m}^{1+}<A_{e,m}^{2+}<\dots<A_{e,m}^{\overline{k}_{i,j,m}^{+}}$ and $0<A_{e,m}^{1-}<A_{e,m}^{2-}<\dots<A_{e,m}^{\overline{k}_{i,j,m}^{-}}$.

\item[Bound constraints]
Bound constraints delineate the domain for each variable. Let $M^+_{i,j,m}$ represent the peak volume of materials that can be excavated from a section $j$ for road $i$ and material $m$, and $M^-_{i,j,m}$ indicate the utmost volume of materials that can be filled into a section $j$ for road $i$ and material $m$. Further, let $\overline{u}_{i,j}$ and $\underline{u}_{i,j}$ signify the uppermost and lowermost offsets for the $j^{th}$ section of the $i^{th}$ road. In a similar vein, $\overline{z}_{e}$ and $\underline{z}_{e}$ stand for the maximum and minimum offset boundaries for intersection $e$. Consequently, the boundary constraints are defined as for all $i \in \mathcal{R}, j \in \mathcal{S}_i, m \in \mathcal{K}, h \in \mathcal{H}$
\begin{gather*}
    0 \leq \var{f^{u+,h}_{i,j,m}} \leq M^{+,h}_{i,j,m},\qquad
    0 \leq \var{f^{u-,h}_{i,j,m}} \leq M^{+,h}_{i,j,m},\qquad
    0 \leq \var{f^{l+,h}_{i,j,m}} \leq M^{+,h}_{i,j,m},\qquad
    0 \leq \var{f^{l-,h}_{i,j,m}} \leq M^{+,h}_{i,j,m},\\
    {M^{+}_{i,j,m} \geq \var{\mathcal{V}^{+,h}_{i,j,m}}} \geq 0, \qquad 
    {M^{-}_{i,j,m} \geq \var{\mathcal{V}^{-,h}_{i,j,m}}} \geq 0; \qquad
    \var{f^{t+,h}_{i,j,m}} \geq 0,\qquad 
    \var{f^{t-,h}_{i,j,m}} \geq 0;
\end{gather*}
for all $e \in \mathcal{I}, m \in \mathcal{K}, h \in \mathcal{H}$,
\begin{gather*}
    {0 \leq \var{f^{u,h}_{e,m}}} \leq M^{+,h}_e,\qquad
    {0 \leq \var{f^{u,h}_{e,m}}} \leq M^{-,h}_{e,m},\qquad
    {M^{+}_{e,m} \geq \var{\mathcal{U}^{+,h}_{e,m}}} \geq 0, \qquad
    {M^{-}_{e,m} \geq \var{\mathcal{U}^{-,h}_{e,m}}} \geq 0;    
\end{gather*}
for all $i \in \mathcal{R}, j \in \mathcal{S}_i$,
\(
    \underline{u}_{i,j} \leq \var{u_{i,j}} \leq \overline{u}_{i,j};
\)
for all $e \in \mathcal{I}$,
\(
    \underline{z}_{e} \leq \var{z_{e}} \leq \overline{z}_{e};
\)
for all $i \in \mathcal{R}, k \in \mathcal{B}_i, m \in \mathcal{K}, h \in \mathcal{H},$
\(
    \var{f^{b+, k, h}_{i,\nu(k),m}} \geq 0, \var{f^{b-,k,h}_{i,\nu(k),m}} \geq 0;
\)
for all $i \in \mathcal{R}, l \in \mathcal{W}_i, m \in \mathcal{K}, h \in \mathcal{H},$
\(
    \var{f^{w+, l, h}_{i,\delta(l),m}} \geq 0, \var{f^{w-, l, h}_{i,\delta(l),m}} \geq 0;
\)
for all $i \in \mathcal{R}, j \in \mathcal{S}_i, e \in \mathcal{I}, m \in \mathcal{K}, h \in \mathcal{H}$,
\(
    \var{f^{in, e, h}_{i,j,m}} \geq 0, \var{f^{out, e, h}_{i,j,m}} \geq 0;
\)
and for all $e \in \mathcal{I}, m \in \mathcal{K},$
\(
    \var{k_{i,j,m}^{+}} \in \{1,2,\dots,\overline{k}_{i,j,m}^{+}\}, 
    \var{k_{i,j,m}^{-}} \in \{1,2,\dots,\overline{k}_{i,j,m}^{-}\},
    \var{k_{e,m}^{+}} \in \{1,2,\dots,\overline{k}_{e,m}^{+}\},
    \var{k_{e,m}^{-}} \in \{1,2,\dots,\overline{k}_{e,m}^{-}\}.
\)
\end{description}

The resulting model is a MILP where the integer variables $\var{k_{i,j,m}^{+}}$, $\var{k_{i,j,m}^{-}}$, $\var{k_{e,m}^{+}}$ and $\var{k_{e,m}^{-}}$ make it a nonconvex model. We next build a convex model.

\section{Convex Unified Vertical Alignment (CUVA) model}\label{s:CUVA}

We first recall Momo et al.~\cite{MOMO-23} convex QCQP model that approximates side slopes using angles and a quadratic form. The model is validated on a single material, one haul path, and only one road. It enjoys a polynomial time worst-case complexity due to convexity (vs. exponential for the MILP model). We then extend the model for multi-material handling.

\subsection{Volume Constraints for single material}
For any given cut or fill section, the volume constraints calculate the volumes of cut or fill by factoring in the side slopes. These constraints confirm that the cumulative volume added to or excavated from a section matches the volume difference between the profiles of the road and ground. Similar to UVA, the elevation or the difference in height between the road and ground profile is depicted as $u_{i,j}$ for every $i \in \mathcal{R}$ and $j \in \mathcal{S}$, with $\underline{u}_{i,j}\leq \var{u_{i,j}}\leq\overline{u}_{i,j}$. Here, $\underline{u}_{i,j}$ and $\overline{u}_{i,j}$ specify the height boundaries within which ground data is accessible. We also note that $\underline{u}_{i,j} < 0$ and $\overline{u}_{i,j} > 0$.

\subsubsection{Side-slopes}
In addressing specific cut or fill sections, we design volume constraints aimed at estimating the volumes for cut or fill, factoring in side slopes. These constraints ensure that the total volume added to or removed from a section matches the volume difference between the road and the ground profile. If the cross-sections are rectangular, the area is $A_{i,j}$ for $i \in \mathcal{R}$ and $j \in \mathcal{S}_{i}$, then we have,
\begin{equation}
	\label{eqn:qcqpVolConsRectangle}
	\begin{aligned}
		\var{\mathcal{V}_{i,j}^{-}} - \var{\mathcal{V}_{i,j}^{+}} = A_{i,j}\var{u_{i,j}}.
	\end{aligned}
\end{equation}

A negative height difference $\var{u_{i,j}}$ indicates a cut section and in this case $\var{\mathcal{V}_{i,j}^{-}} = 0$ in \eqref{eqn:qcqpVolConsRectangle}. Similarly, for a fill section, the height difference, $\var{u_{i,j}}$ is positive and $\var{\mathcal{V}_{i,j}^{+}} = 0$ in \eqref{eqn:qcqpVolConsRectangle}. However, when side slopes are included, the cross-sectional shape of the road at a specific section becomes trapezoidal. Consequently, the volume becomes a quadratic function in relation to $\var{u_{i,j}}$, thereby breaking the model's linear structure as shown in \cite{HARE-15}. To preserve the model's linearity, side slopes are approximated using rectangular slabs stacked atop each other in the MILP model.

\begin{figure}
	\centering
	\includegraphics[width=0.5\textwidth,keepaspectratio]{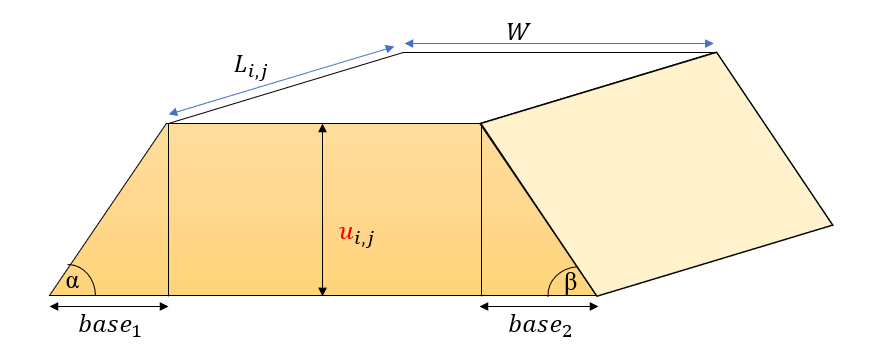}
	\caption{A typical section with side-slopes and trapezoid cross-section in the QCQP model}
	\label{fig:QCQP-UVA-side-slope-1}
\end{figure}

In the QCQP model, side-slopes are not approximated. Instead, the angles of the side slopes are approximated and are in turn used to directly compute the trapezoid's area. The angles are approximated using least-squares for each section $j \in \mathcal{S}_{i}$. By increasing the number of rectangular slabs in the MILP model, the combined area of these slabs approaches the true area of the trapezoid. This relationship serves as the primary link between the two models. Within the QCQP model, volume constraints are established by determining the volume formed by two triangles and a rectangle. Figure~\ref{fig:QCQP-UVA-side-slope-1} provides a visual representation of this idea. In this depiction, $\alpha$ represents the angle of the left side-slope, while $\beta$ signifies the angle of the right side-slope.

The area of the left (resp. right) triangle is given by $\var{u}^2 / (2 \tan \alpha)$ (resp. $\var{u}^2 / (2 \tan \beta)$). So, the area of the cross-section is 
$\var{u} W + 1/2 \var{u}^{2}\Big(1/(\tan \alpha)+1/(\tan \beta)\Big)$, 
where the area of the rectangle is $\var{u}W$. Noting $A_{i,j} = WL_{i,j}$, for a section $j$ for road $i$, the volume $\var{\mathcal{V}_{i,j}}$ of a section becomes
\begin{equation}
	\label{eqn:qcqpVolConsVol-1}
	\var{\mathcal{V}_{i,j}} = A_{i,j}\var{u_{i,j}} +\frac{1}{2}\var{u_{i,j}}^{2}\Big(\frac{1}{\tan \alpha_{i,j}}+\frac{1}{\tan \beta_{i,j}}\Big)L_{i,j}.
\end{equation}
Hence, \eqref{eqn:qcqpVolConsVol-1} replaces \eqref{eqn:qcqpVolConsRectangle} giving
\[\var{{\mathcal{V}_{i,j}}^{-}} - \var{{\mathcal{V}_{i,j}}^{+}} = A_{i,j}\var{u_{i,j}} +\frac{1}{2}\var{u_{i,j}}^{2}\Big(\frac{1}{\tan \alpha_{i,j}}+\frac{1}{\tan \beta_{i,j}}\Big)L_{i,j}.\]

For cut sections, the road profile lies beneath the ground level, whereas for fill sections, the road profile is elevated above the ground. Since the volume is always non-negative, a viable convex region can be established using inequalities. The volume constraint for a cutting section satisfies
\begin{equation}
	\var{\mathcal{V}_{i,j}^{+}} \geq -A_{i,j}\var{u_{i,j}} +  \frac{1}{2}L_{i,j}\Big(\frac{1}{\tan \big(\alpha_{i,j,cut}\big) }+\frac{1}{\tan \big(\beta_{i,j,cut}\big)}\Big)\var{u_{i,j}}^{2},\ \text{for}\ \underline{u}_{i,j} \leq  \var{u_{i,j}} \leq 0,
\end{equation}
whereas for a fill section
\begin{equation}
	\var{\mathcal{V}_{i,j}^{-}} \geq A_{i,j}\var{u_{i,j}} + \frac{1}{2}L_{i,j}\Big(\frac{1}{\tan \big(\alpha_{i,j,fill}\big) }+\frac{1}{\tan \big(\beta_{i,j,fill}\big)}\Big)\var{u_{i,j}}^{2},\ \text{for}\  0 \leq \var{u_{i,j}} \leq \overline{u}_{i,j}.
\end{equation}
The QCQP model's volume constraints are formulated as quadratic functions, distinguishing them from the volume constraints found in the MILP model. Momo et al.~\cite{MOMO-23} confirm the convex nature of these volume constraints. Given that the remaining constraints and the objective function are linear, this ensures the entire model is convex. The authors approximate values for angles $\alpha$ and $\beta$ for every section $j \in \mathcal{S}_{i}$ using least-square. 

Our focus is now on extending these volume constraints to encompass multiple materials, all while preserving the model's convex property.

\subsection{Volume Constraints for multi-material}
\begin{figure}
	\centering
	\includegraphics[width=0.5\textwidth,keepaspectratio]{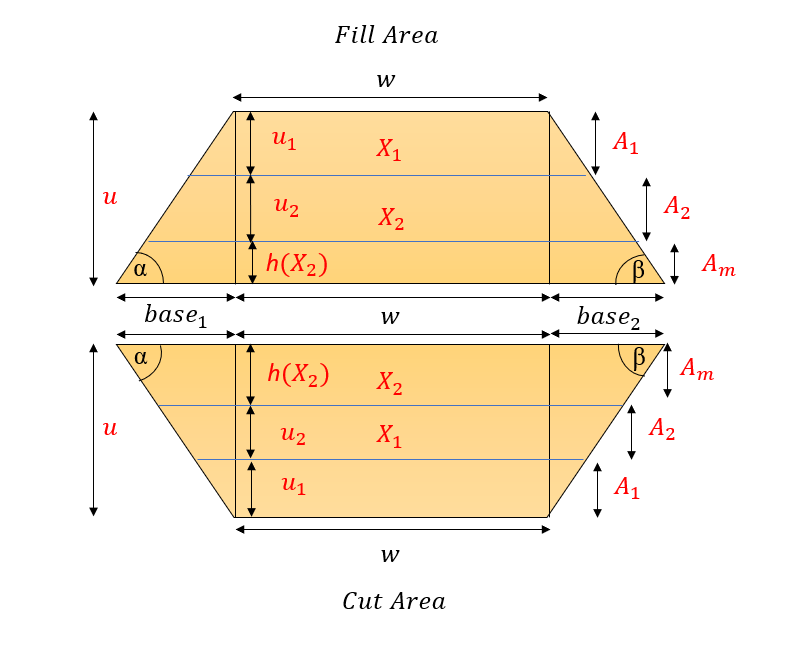}
	\caption{A cross-section of a cut section and fill section with side-slopes when the number of materials equals 3.}
	\label{fig:QCQP-UVA-linear-interpol-2}
\end{figure}

First, we show that extending Momo et al.~\cite{MOMO-23} model to multimaterial breaks convexity. Consider Figure~\ref{fig:QCQP-UVA-linear-interpol-2} where we note $\var{u}$ the height of trapezoid, $A_1$ the area of Material 1 with height $\var{u_1}$ and trapezoid base $\var{X_1}$. Elementary geometry provides the volume for each material as the following constraints
\begin{align}	
	\var{V_m^+} &\geq \frac{1}{2C}(w^2) + \var{u}w + \frac{1}{2}\var{u}^2 - \frac{1}{2C}\var{X_{m-1}}^2,
	\label{eqn:linear-interpol-vol-constraint-8} \\
	\var{V_i^+} &\geq \frac{1}{2C}(\var{X_{i}}^2 - \var{X_{i-1}}^2),
	\label{eqn:linear-interpol-vol-constraint-9} \\	
	\var{V_1^+} &\geq \frac{1}{2C}(\var{X_1}^2 - w^2). \label{eqn:linear-interpol-vol-constraint-10}
\end{align}
It is clear that constraints \eqref{eqn:linear-interpol-vol-constraint-8} and \eqref{eqn:linear-interpol-vol-constraint-9} are not convex.

\subsection{A different approach to formulating volume constraints}

Instead of approximating the angles $\alpha$ and $\beta$ for the trapezoid for a section $j$ as shown earlier, it is possible to approximate the volume directly for a given section $j$ using least-squares. For a single material use case, and when applying linear least squares, the volume constraint for a cut section and a fill section is expressed
\begin{equation}
	\var{\mathcal{V}_{i,j}^{+}} \geq \chi_{i,j,1}\var{u_{i,j}} +  \chi_{i,j,2},\ \text{for}\ \underline{u}_{i,j} \leq  \var{u_{i,j}} \leq 0,
\end{equation}
\begin{equation}
	\var{\mathcal{V}_{i,j}^{-}} \geq \chi_{i,j,1}\var{u_{i,j}} +  \chi_{i,j,2},\ \text{for}\ 0 \leq  \var{u_{i,j}} \leq \overline{u}_{i,j},
\end{equation}
for all $i \in \mathcal{R}$, and $j \in \mathcal{S}_i.$ Similarly, when applying quadratic least squares, the volume constraint for a cut and fill section is expressed
\begin{equation}
	\var{\mathcal{V}_{i,j}^{+}} \geq \chi_{i,j,1}\var{u_{i,j}^2} +  \chi_{i,j,2}\var{u_{i,j}} + \chi_{i,3}, \ \text{for}\ \underline{u}_{i,j} \leq  \var{u_{i,j}} \leq 0,
\end{equation}
\begin{equation}
	\var{\mathcal{V}_{i,j}^{-}} \geq \chi_{i,j,1}\var{u_{i,j}^2} +  \chi_{i,j,2}\var{u_{i,j}} + \chi_{i,j,3}, \ \text{for}\ 0 \leq  \var{u_{i,j}} \leq \overline{u}_{i,j},
\end{equation}
where  $\forall i \in \mathcal{R}$, and $\forall j \in \mathcal{S}_i.$ Here, $\chi_{i,1}$ and $\chi_{i,2}$ are the coefficients of the linear least squares algorithm for section $i$ and $\chi_{i,1}$, $\chi_{i,2}$, and $\chi_{i,3}$ are the coefficients of the quadratic least squares algorithm for section $i$.

Now we incorporate it into our CUVA model, extending the volume constraints to multiple materials, we have the constraints as follows when using least squares with a linear model fit for cut and fill sections 
\begin{equation}
	\var{\mathcal{V}_{i,j,m}^{+}} \geq \chi_{i,j,m,1}\var{u_{i,j,m}} +  \chi_{i,j,m,2},\ \text{for}\ \underline{u}_{i,j,m} \leq  \var{u_{i,j,m}} \leq 0,
\end{equation}
\begin{equation}
	\var{\mathcal{V}_{i,j,m}^{-}} \geq \chi_{i,j,m,1}\var{u_{i,j,m}} +  \chi_{i,j,m,2},\ \text{for}\ 0 \leq  \var{u_{i,j,m}} \leq \overline{u}_{i,j,m},
\end{equation}
for all $i \in \mathcal{R}$, $j \in \mathcal{S}_i$, and $m \in \mathcal{K}$. Similarly, when using least squares with a quadratic model fit, for cut and fill sections, we have
\begin{equation}
	\var{\mathcal{V}_{i,j,m}^{+}} \geq \chi_{i,j,m,1}\var{u_{i,j,m}}^2 +  \chi_{i,j,m,2}\var{u}_{i,j,m} + \chi_{i,j,m,3}, \ \text{for}\ \underline{u}_{i,j,m} \leq  \var{u_{i,j,m}} \leq 0,
\end{equation}
\begin{equation}
	\var{\mathcal{V}_{i,j,m}^{-}} \geq \chi_{i,j,m,1}\var{u_{i,j,m}}^2 +  \chi_{i,j,m,2}\var{u}_{i,j,m} + \chi_{i,j,m,3}, \ \text{for}\ 0 \leq  \var{u_{i,j,m}} \leq \overline{u}_{i,j,m},
\end{equation}
for all  $i \in \mathcal{R}$, $j \in \mathcal{S}_i$, $m \in \mathcal{K}$.
In a similar manner, the constraints for the intersections for a cut and fill site respectively, when applying a linear model are
\begin{equation}
	\var{\mathcal{U}_{e,m}^{+}} \geq \chi_{e,m,1}\var{z_{e,m}} +  \chi_{e,m,2},\ \text{for}\ \underline{z}_{e,m} \leq  \var{z_{e,m}} \leq 0,
\end{equation}
\begin{equation}
	\var{\mathcal{U}_{e,m}^{-}} \geq \chi_{e,m,1}\var{z_{e,m}} +  \chi_{e,m,2},\ \text{for}\ 0 \leq  \var{z_{e,m}} \leq \overline{z}_{e,m},
\end{equation}
for all  $e \in \mathcal{I}$, $m \in \mathcal{K}$; and the constraints for the intersections when applying a quadratic model are
\begin{equation}
	\var{\mathcal{U}_{e,m}^{+}} \geq \chi_{e,m,1}\var{z_{e,m}}^2 +  \chi_{e,m,2}\var{z}_{e,m} + \chi_{e,m,3}, \ \text{for}\ \underline{z}_{e,m} \leq  \var{z_{e,m}} \leq 0,
\end{equation}
\begin{equation}
	\var{\mathcal{U}_{e,m}^{-}} \geq \chi_{e,m,1}\var{z_{e,m}}^2 +  \chi_{e,m,2}\var{z}_{e,m} + \chi_{e,m,3}, \ \text{for}\ 0 \leq  \var{z_{e,m}} \leq \overline{z}_{e,m},
\end{equation}
for all $e \in \mathcal{I}$, and $m \in \mathcal{K}$.

\subsection{Convexity}
\begin{figure}[ht]
	\centering
	\includegraphics[width=0.6\textwidth,keepaspectratio]{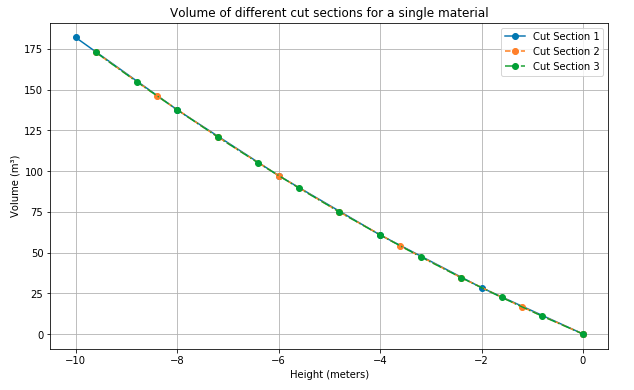}
	\caption{Volume of different cut sections for a single material.}
	\label{fig:vol-cut-sec-conv-det}
\end{figure}

\begin{figure}[ht]
	\centering
	\includegraphics[width=0.6\textwidth,keepaspectratio]{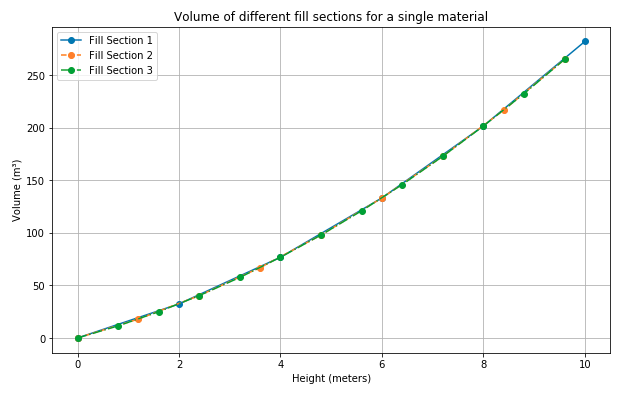}
	\caption{Volume of different fill sections for a single material.}
	\label{fig:vol-fill-sec-conv-det}
\end{figure}

Our objective function and all our constraints except for the volume constraints are linear. Our volume constraints have linear equations when using linear models and quadratic equations when using quadratic models. So, to confirm convexity, only the quadratic models need to be examined.
Figure \ref{fig:vol-cut-sec-conv-det} and \ref{fig:vol-fill-sec-conv-det} illustrate the volume function for different cut and fill sections respectively, for a single material. 

\begin{figure}
	\centering
	\includegraphics[width=0.6\textwidth,keepaspectratio]{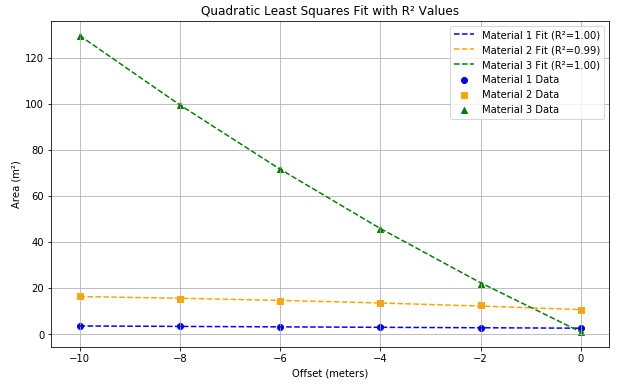}
	\caption{Quadratic least squares approximation for a cut section for multiple materials.}
	\label{fig:quadratic-least-squares-cut-mul}
\end{figure}

\begin{figure}
	\centering
	\includegraphics[width=0.6\textwidth,keepaspectratio]{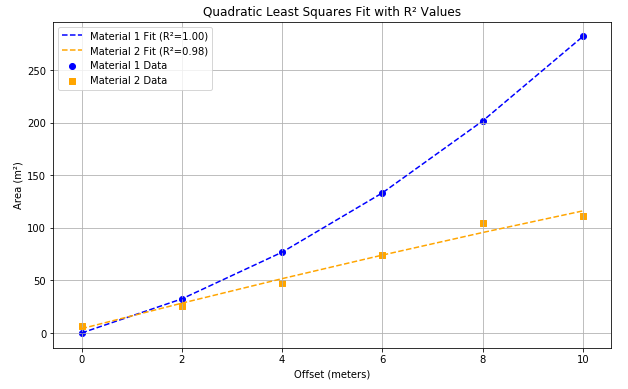}
	\caption{Quadratic least squares approximation for a fill section for multiple materials.}
	\label{fig:quadratic-least-squares-fill-mul}
\end{figure}

The volume constraint for a cut section is 
$$\var{\mathcal{V}^{+}} \geq \chi_{1}\var{u}^2 +  \chi_{2}\var{u} + \chi_{3}.$$
Similarly, for the intersection volume constraints, the volume constraint for a cut section is 
$$\var{\mathcal{U}^{+}} \geq \chi_{1}\var{z}^2 +  \chi_{2}\var{z} + \chi_{3}.$$
Figures \ref{fig:quadratic-least-squares-cut-mul} and \ref{fig:quadratic-least-squares-fill-mul} show the quadratic least square approximation for a cut and fill section respectively, for multiple materials.

In our implementation, we follow \cite{MOMO-23} approach and use $R^2$ to select the between linear and quadratic models. This approach maintains the convexity of the approximation for all road sections within our dataset. In contrast, relying solely on quadratic models results in  breaching the convexity of some constraints for three roads in our dataset. 
A quadratic least-square fit is typically effective since volumes often have a trapezoidal shape, leading to a quadratic increase of volume with respect to the offset.

\section{Numerical Experiments}\label{s:val}
In this section, we assess the validity of our CUVA model by juxtaposing its results against those from the QCQP-QNF, MR-QNF, and MH-QNF, using data pertinent to each model. Additionally, to further ensure the model's robustness, we conducted comprehensive unit tests, including all those referenced in \cite{MOMO-23,AYMAN-23,BEIRANVAND-17}.

All experimental tests were conducted on a consistent setup, equipped with an AMD Ryzen 7 3700X 8-Core 3.60 GHz Processor and 80 GB RAM. Both our UVA model and the MR-QNF model were developed using MATLAB R2022a 64-bit edition \cite{MATLAB-21} on a Windows 11 platform. The MH-QNF model was crafted in C++ and compiled with the MSVC x64 compiler in Visual Studio 2022, also on Windows 11. For optimization variable and parameter modeling in MATLAB, we utilized the YALMIP \cite{LOEFBERG-04} optimization toolbox (version R20230609) in conjunction with the academic release of the IBM ILOG CPLEX solver V12.10.0 \cite{CPLEX-09}. We opted for this earlier version of CPLEX as the later releases no longer provide support for running on MATLAB. Additionally, we experimented with the latest version of Gurobi \cite{ACHTERBERG-13} and noticed that CPLEX consistently delivered equal or superior results for our dataset.

For our CUVA model, we leverage the input data format from MR-QNF \cite{AYMAN-23}. This data format is applied to both our single-road and multi-road scenarios, and we have extended it to accommodate multiple materials.

A crucial element in our experiments is the cost components incorporated within the CUVA model. To ensure a balanced comparison of our model's performance, we align these components with those utilized in the MR-QNF and MH-QNF models. 

\subsection{QCQP-QNF vs. CUVA}

We first compare our technique of approximating the volume of materials in the CUVA model against the QCQP-QNF model~\cite{MOMO-23}, for a single material use case. In contrast to the CUVA model, the QCQP-QNF model cannot deal with road networks, multiple materials, or multiple hauling paths. 

For this comparison, we take the mean absolute percentage error (MAPE) as well as the root mean squared error (RMSE). For computing the MAPE for a road, we first take the absolute difference between the actual area and the predicted area and compute the percentage error, averaged over each offset in our data, and finally averaged over all the road sections. This is computed separately for cut and fill sections of a road. The MAPE for a road for all cut sections is given by
\begin{equation}
	\label{eqn:MAPE-cut}
	MAPE^+ = \frac{1}{n_{S}}\sum_{i\in\mathcal{S}}\frac{1}{n_{u_{i}}}\sum_{j\in u_{i}}\frac{|P^{+}_{j}-\hat{y_{j}^+}|}{P^{+}_{j}} 100\%.
\end{equation}
Similarly, the MAPE for all fill sections of a road is given by
\begin{equation}
	\label{eqn:MAPE-fill}
	MAPE^- = \frac{1}{n_{S}}\sum_{i\in\mathcal{S}}\frac{1}{n_{u_{i}}}\sum_{j\in u_{i}}\frac{|P^{-}_{j}-\hat{y_{j}^-}|}{P^{-}_{j}} 100\%.
\end{equation}

Here, $\hat{y^+}$ and $\hat{y^-}$ are the predicted area computed by the least squares model for a given offset in a cut and fill section respectively. The values of $\hat{P^{+}}$ and $\hat{P^{-}}$ are the actual area for a given offset in a cut and fill section respectively, as extracted from RoadEng software (in the form of CSV files) provided by Softree Technical Systems Inc\footnote{\url{softree.com}}. We also take the root mean squared error (RMSE) by calculating the squared error between the actual and the predicted area averaged over all the offsets for a given section and then compute its square root. We then average this value over all the sections of a road. The RMSE for a road for all cut sections is given by
\begin{equation}
	\label{eqn:RMSE-cut}
	RMSE^+ = \frac{1}{n_{S}}\sum_{i\in\mathcal{S}}\sqrt{\frac{1}{n_{u_{i}}}\sum_{j\in u_{i}}(P^{+}_{j}-\hat{y_{j}^+})^2}.
\end{equation}
Similarly, the RMSE for a road for all fill sections is given by
\begin{equation}
	\label{eqn:RMSE-fill}
	RMSE^- = \frac{1}{n_{S}}\sum_{i\in\mathcal{S}}\sqrt{\frac{1}{n_{u_{i}}}\sum_{j\in u_{i}}(P^{-}_{j}-\hat{y_{j}^-})^2}.
\end{equation}

MAPE expresses errors as a percentage of the actual values, providing a sense of the size of the errors in a relative term. We compute the MAPE and RMSE for 52 roads found in \cite{MOMO-23} for cut and fill sections separately, for both the QCQP-QNF model and our CUVA model. We set the significance threshold for our experiments at $0.05$. This is summarized in Table \ref{tab:qcqp-qnf-cut-error} and Table \ref{tab:qcqp-qnf-fill-error} for cut and fill sections respectively (in Appendix \ref{appendix:qcqp-old-new-comp}). We truncate the values up to two decimal places for better readability.

\begin{figure}
	\centering
	\includegraphics[width=0.8\textwidth]{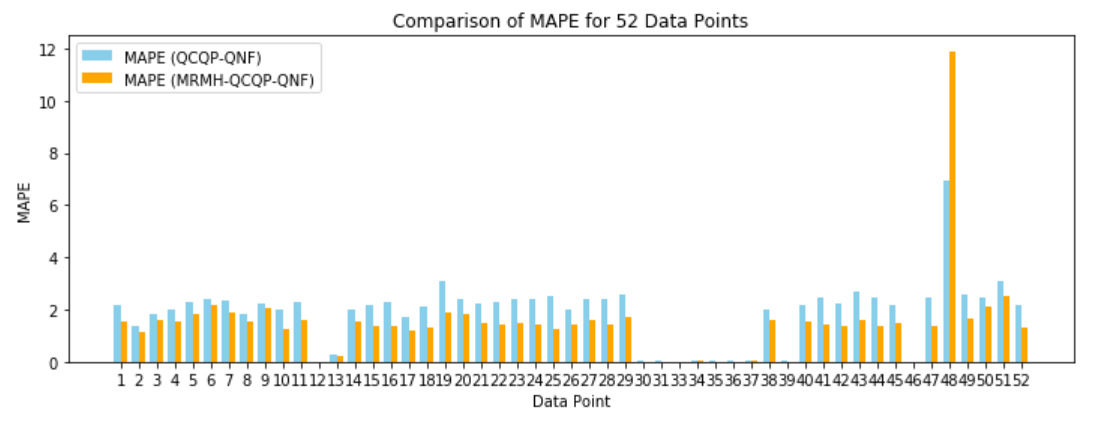}
	\caption{Comparison of MAPE between QCQP-QNF and CUVA(ours) across all FILL sections for 52 road dataset - excluding outliers.}
	\label{fig:MAPE-fill}
\end{figure}

The comparison of the MAPE and RMSE between QCQP-QNF and CUVA shows that the later has a lower MAPE. An  illustration is provided in Figure~\ref{fig:MAPE-fill}.
For enhanced clarity, outliers (such as road $32$) have been excluded from the analysis. 


%

In these tables and figures, we observe that for all the cut sections, the MAPE and the RMSE are significantly lower for our CUVA model for all the roads. For the fill sections, we observe a similar trend. However, for roads $32$ and $48$, the MAPE was observed to be better for the QCQP-QNF model, and for roads $9$, $50$, and $51$, the RMSE was observed to be better for the QCQP-QNF model. On further investigation, we find that for the CUVA model, for roads $48$, $51$, and $52$, we have $585$, $51$, and $165$ sections respectively where the quadratic approximation results in poor $R^2$ values or higher residual errors. For these cases, the model makes use of the linear approximation instead since it leads to better $R^2$ values. 

\begin{figure}
	\centering
	\begin{subfigure}[b]{0.45\textwidth}
		\centering
		\includegraphics[width=\textwidth]{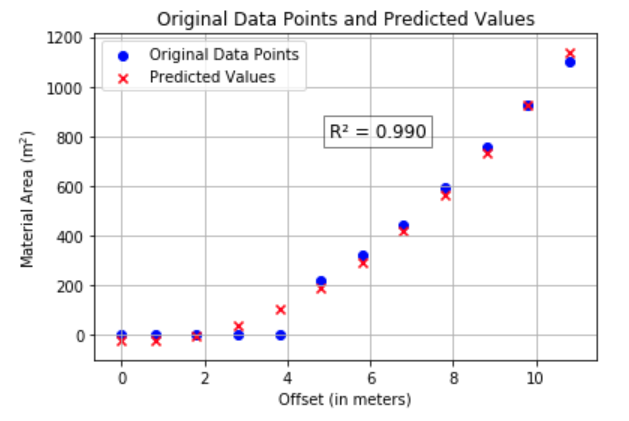}
		\caption{Quadratic approximation of fill material for Road $32$ Section $120$.}
		\label{fig:approx-analysis-1}		
	\end{subfigure}
	\hfill
	\begin{subfigure}[b]{0.45\textwidth}
		\centering
		\includegraphics[width=\textwidth]{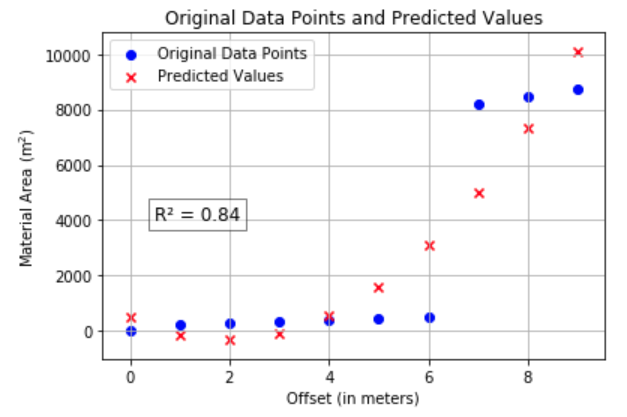}
		\caption{Quadratic approximation of cut material for Road $48$ Section $1912$.}
		\label{fig:approx-analysis-2}
	\end{subfigure}
	\hfill
	\begin{subfigure}[b]{0.45\textwidth}
		\centering
		\includegraphics[width=\textwidth]{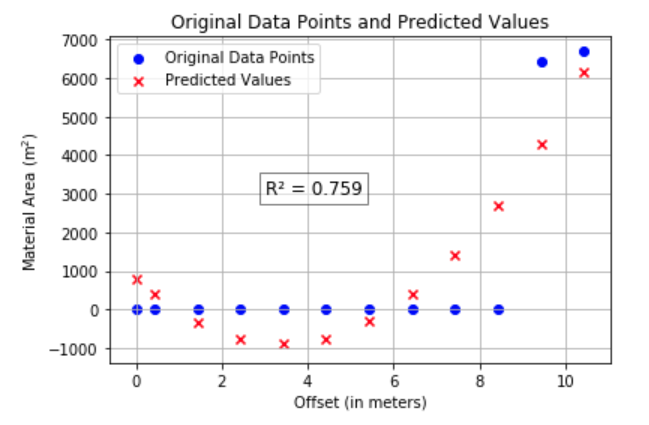}
		\caption{Quadratic approximation of fill material for Road $48$ Section $1912$.}
		\label{fig:approx-analysis-3}
	\end{subfigure}
	\caption{Quadratic approximation of material, and the corresponding $R^2$ value.}
	\label{fig:bplq}
\end{figure}

In contrast, for the QCQP-QNF model, for road $48$, $51$, and $52$, we have $721$, $786$, and $535$ sections respectively for which the quadratic approximation leads to lower $R^2$ values. Additionally, for $12$ other roads (6, 7, 9, 11, 14, 29, 32, 38, 42, 43, 47, 50), the QCQP-QNF model's quadratic approximation leads to lower $R^2$ values. However, this is not the case for the CUVA model. Therefore, defaulting to a linear approximation for cases when the quadratic approximation leads to a higher error, fails to explain the higher error for roads $9$ and $32$. We investigate the fill section from road $32$ with the maximum error as well as the cut and fill sections from road $48$ with the maximum error. The quadratic approximations of the area of material for these sections are plotted in Figure \ref{fig:approx-analysis-1}, Figure \ref{fig:approx-analysis-2}, and Figure \ref{fig:approx-analysis-3}.

For road $32$ the $R^2$ value in Figure~\ref{fig:approx-analysis-1} is high irrespective of higher MAPE for the overall road. We find that for road $48$, the original data distribution is not quadratic. According to our industrial partner, this kind of distribution of the material profile is unusual, and is caused by the cross-section calculation. This can be an avenue for future research.




\subsection{MH-QNF vs. CUVA}

In this analysis, we evaluate the performance of both the MH-QNF model and the CUVA model using the road data from \cite{BEIRANVAND-17}. In contrast to the CUVA model, the MH-QNF model cannot deal with road networks. This data contains $23$ roads after omitting roads with blocks or obstacles. Additionally, we examine the execution times required to solve each road data set (solver time), calculating the average across five runs. For this time comparison, we use the ratio of the base model's duration to that of the new model. These comparisons are showcased in Table \ref{tab:MH-QNF-optimal-cost-comparison} in Appendix \ref{appendix:qcqp-milp-comp}. 


Remarkably, for road 7, while the MH-QNF model yields an infeasible solution, the CUVA model successfully provides an optimal solution. We also note that for all roads, the MH-QNF model was marginally faster, however, we used C++ for running the MH-QNF model as compared to MATLAB for running the CUVA model. The speedup in solve time for the MH-QNF model decreased by increasing the size of the road. In \cite{MOMO-23}, the authors showed that the QCQP-QNF model was faster as compared to the MH-QNF model. 

%

\subsection{MR-QNF vs. CUVA}

We proceed to compare the results of our CUVA model with that of the MR-QNF model using simpler road network data that we call the academic data, as well as real-world road network data provided in \cite{AYMAN-23}. In contrast to the CUVA model, the MR-QNF model cannot deal with multiple materials or multiple hauling paths. This comparison is shown in Table \ref{tab:road_network_QCQP_acad_optimal_cost} for academic data, and Table \ref{tab:road_network_QCQP_real_world_optimal_cost} for real-world data. 
The optimal cost in these tables is rounded to the nearest integer.

In the case of academic road network data, which typically features smaller networks compared to actual road systems, no significant difference was noticed for the solve time between the CUVA and the MR-QNF model. 

Turning to real-world road network data, we utilize two road networks consisting of logging roads. The first network comprises twelve roads and five intersections, spanning approximately 6 kilometers. Here, although the MR-QNF model produces a lower optimal cost, it requires double the solver time of the CUVA model. It was noted in \cite{MOMO-23} that the MILP model results in better optimal costs for roads with a larger number of slabs. The second real-world example encompasses a larger network with 32 roads and 16 intersections, extending over roughly 66 kilometers. In this scenario, the MR-QNF model fails to compute a solution after an extended duration (timed out/out of memory after 30 minutes), whereas the CUVA model efficiently delivers an optimal cost of $174,610$ in about 9 seconds. We note that this optimal cost is significantly lower than the one reported in \cite{AYMAN-23} using the one-at-a-time method. 

		\begin{table}
			\centering
			\caption{Comparison of solve time for MR-QNF and CUVA model for academic road network data for single material.}
			\label{tab:road_network_QCQP_acad_optimal_cost}
            \begin{tabular}{@{}crrrrrrr@{}}
				\toprule
				\thead{Road\\Network \\ \#} & \thead{Roads} & \thead{Stations} & \thead{Intersections} & \thead{Slabs} 
    & \thead{Time \\in \\seconds \\(MR-QNF)} & \thead{Time \\in \\seconds \\(MR-QNF)} & \thead{Speedup} \\ \midrule
				1              & 3     & 9        & 1             & 5     
    & 0.04 & 0.05 & 0.08\\
				2              & 3     & 9        & 1             & 5     
    & 0.05 & 0.05 & 1\\
				3              & 4     & 12        & 1             & 5     
    & 0.05 & 0.05 & 1\\
				\bottomrule
			\end{tabular}
		\end{table}	

		\begin{table}
            \centering
			\caption{Comparison of solve time for MR-QNF and CUVA model for real-world road network data for single material.}
			\label{tab:road_network_QCQP_real_world_optimal_cost}
			\begin{tabular}{@{}crrrrrrr@{}}
				\toprule
				\thead{Road\\Network \\ \#} & \thead{Roads} & \thead{Stations} & \thead{Intersections} & \thead{Slabs} 
    & \thead{Time \\in \\seconds \\(MR-QNF)} & \thead{Time \\in \\seconds \\(CUVA)} & \thead{Speedup} \\ \midrule
				1              & 12    & 327      & 5             & 42    
    & 0.6 & 0.2 & 2.7 \\
				2              & 32    & 1,508     & 16            & 22    
    & out-of-memory & 8.9 & -       \\
				\bottomrule
			\end{tabular}
		\end{table}

\subsection{UVA vs. CUVA}

We now proceed to compare the results of our CUVA model and UVA model using multi-material road network data, which is derived from the road network data provided in \cite{AYMAN-23} and extended to multiple layers of materials synthetically. In this case, the data contains three different types of materials. This comparison is shown in Table \ref{tab:road_network_acad_optimal_cost_final} for academic data, and Table \ref{tab:road_network_real_world_optimal_cost_final} for real-world data. 
The optimal cost in these tables is rounded to the nearest integer.

%


Examining real-world road networks with multiple materials, the initial example includes a network of twelve roads and five intersections, spanning a total length of 6 kilometers, while the real-world network features an extensive network with 32 roads and 16 intersections, spanning an approximate length of 66 kilometers. For the first real-world road, the UVA model generates a significantly lower ($-7.21\%$) optimal cost as compared to the CUVA model, however, the UVA model takes 16 times more computation time to solve the problem. For the second real-world road data, the UVA model fails to give an optimal solution after a prolonged computation period and runs out of memory. We timed out the computation at 30 minutes. The CUVA model on the other hand gives us an optimal solution in a reasonable time. We can conclude that for very large networks, we can rely on the CUVA model to give us an optimal solution.

		\begin{table}
			\centering
			\caption{Comparison of solve time for UVA and CUVA model for academic road network data for multiple materials.}
			\label{tab:road_network_acad_optimal_cost_final}
			\begin{tabular}{@{}crrrrrrr@{}}
				\toprule
				\thead{Road\\Network \\ \#} & \thead{Roads} & \thead{Stations} & \thead{Intersections} & \thead{Slabs} 
    & \thead{Time \\in \\seconds \\(UVA)} & \thead{Time \\in \\seconds \\(CUVA} & \thead{Speedup} \\ \midrule
				1              & 3     & 9        & 1             & 5     
    & 0.05 & 0.05 & 1 \\
				2              & 3     & 9        & 1             & 5     
    & 0.05 & 0.05 & 1 \\
				3              & 4     & 12        & 1             & 5     
    & 0.05 & 0.05 & 1 \\
				\bottomrule
			\end{tabular}
		\end{table}

        \begin{table}
            \centering
			\caption{Comparison of solve time for UVA and CUVA model for real-world road network data for multiple materials.}
			\label{tab:road_network_real_world_optimal_cost_final}
			\begin{tabular}{@{}crrrrcrr@{}}
				\toprule
				\thead{Road\\Network \\ \#} & \thead{Roads} & \thead{Stations} & \thead{Intersections} & \thead{Slabs} 
    & \thead{Time \\in \\seconds \\(UVA)} & \thead{Time \\in \\seconds \\(CUVA)} & \thead{Speedup} \\ \midrule
				1              & 12    & 327      & 5             & 42    
    & 18 & \textbf{1.1} & \textbf{16.4} \\
				2              & 32    & 1,508     & 16            & 22    
    & out-of-memory & \textbf{61} & -       \\
				\bottomrule
			\end{tabular}
		\end{table}

\section{Conclusion}\label{s:conc}

In the present study, we have expanded upon the QNF model as delineated in \cite{HARE-15} by integrating the concepts of multiple-hauling paths \cite{BEIRANVAND-17}, multiple materials \cite{HOSSAIN-13} and multiple roads in a road network \cite{AYMAN-23}, while also considering the impact of side slopes. In developing our UVA model,  we observed that the model is prone to encountering memory overflows when applied to extensive road networks, akin to the instances described in \cite{AYMAN-23}. This limitation arises from the intrinsic characteristics of the mixed-integer linear programming (MILP) approach, which has a worst-case exponential time complexity.

To address the computational challenges associated with the MILP framework, we explored the utilization of quadratic constraints as proposed by \cite{MOMO-23} to obtain a convex model. Our investigations revealed that the model's convexity of~\cite{MOMO-23} is compromised when attempts are made to extend it to scenarios involving multiple materials. Consequently, we introduce an alternative approach to formulate the volume constraints for a single material scenario that maintains the convex nature of the model. This method was successfully adapted to the multi-material context, ensuring the preservation of convexity. Convex QCQP models exhibit a worst-case time complexity of polynomial time, positioning them as a faster alternative to MILP models. 

The CUVA model can process large road networks more swiftly than the UVA model. Additionally, for even more extensive road networks, the CUVA model manages to compute solutions within a reasonable time, whereas the UVA model fails to solve such problems.

\section*{Acknowledgments}
The research was supported by the Natural Sciences and Engineering Research Council of Canada under Collaborative Research and Development grant CRDPJ 479316 - 15, and by Innovate BC under Ignite grant IGNITE-2021-RND11-277-UBCO-Lucet-Softree. Both grants were supported by Softree Technical Systems Inc. who was instrumental in providing data (ground maps), software license for RoadEng, and technical expertise.

\section*{Data availability statement}
Due to the nature of the research, commercial supporting data is not available.

\printbibliography

@Book{JHA-06,
  author    = {Jha, M. K. and Schonfeld, P. and Jong, J. C.},
  publisher = {WIT press},
  title     = {Intelligent road design},
  year      = {2006},
  volume    = {19},
}

@Article{HARE-11,
  author    = {Hare, W. L. and Koch, V. R. and Lucet, Y.},
  journal   = {European Journal of Operational Research},
  title     = {Models and algorithms to improve earthwork operations in road design using mixed integer linear programming},
  year      = {2011},
  number    = {2},
  pages     = {470--480},
  volume    = {215},
  publisher = {Elsevier},
}

@Article{MONDAL-15,
  author    = {Mondal, S. and Lucet, Y. and Hare, W.},
  journal   = {Computers \& Operations Research},
  title     = {Optimizing horizontal alignment of roads in a specified corridor},
  year      = {2015},
  pages     = {130--138},
  volume    = {64},
  publisher = {Elsevier},
}

@Article{EASA-88,
  author    = {Easa, S. M.},
  journal   = {Transportation Research Part A: General},
  title     = {Selection of roadway grades that minimize earthwork cost using linear programming},
  year      = {1988},
  number    = {2},
  pages     = {121--136},
  volume    = {22},
  publisher = {Elsevier},
}

@Article{MOREB-96,
  author    = {Moreb, A. A.},
  journal   = {European Journal of Operational Research},
  title     = {Linear programming model for finding optimal roadway grades that minimize earthwork cost},
  year      = {1996},
  number    = {1},
  pages     = {148--154},
  volume    = {93},
  publisher = {Elsevier},
}

@Article{MOREB-04,
  author    = {Moreb, A. A. and Aljohani, M. S.},
  journal   = {Journal of Systems Science and Systems Engineering},
  title     = {Quadratic representation for roadway profile that minimizes earthwork cost},
  year      = {2004},
  number    = {2},
  pages     = {245--252},
  volume    = {13},
  publisher = {Springer},
}

@Article{KOCH-10,
  author    = {Koch, V. R. and Lucet, Y.},
  journal   = {Journal of Industrial \& Management Optimization},
  title     = {A note on: Spline technique for modeling roadway profile to minimize earthwork cost},
  year      = {2010},
  number    = {2},
  volume    = {6},
  publisher = {American Institute of Mathematical Sciences},
}

@Article{HARE-15,
  author    = {Hare, W. and Lucet, Y. and Rahman, F.},
  journal   = {European Journal of Operational Research},
  title     = {A mixed-integer linear programming model to optimize the vertical alignment considering blocks and side-slopes in road construction},
  year      = {2015},
  number    = {3},
  pages     = {631--641},
  volume    = {241},
  publisher = {Elsevier},
}

@Article{BEIRANVAND-17,
  author    = {Beiranvand, V. and Hare, W. and Lucet, Y. and Hossain, S.},
  journal   = {Engineering Optimization},
  title     = {Multi-haul quasi network flow model for vertical alignment optimization},
  year      = {2017},
  number    = {10},
  pages     = {1777--1795},
  volume    = {49},
  publisher = {Taylor \& Francis},
}

@Article{LEE-01,
  author    = {Lee, Y. and Cheng, J.},
  journal   = {Journal of Transportation Engineering},
  title     = {Optimizing highway grades to minimize cost and maintain traffic speed},
  year      = {2001},
  number    = {4},
  pages     = {303--310},
  volume    = {127},
  publisher = {American Society of Civil Engineers},
}

@Article{FWA-02,
  author    = {Fwa, T. F. and Chan, W. T. and Sim, Y. P.},
  journal   = {Journal of Transportation Engineering},
  title     = {Optimal vertical alignment analysis for highway design},
  year      = {2002},
  number    = {5},
  pages     = {395--402},
  volume    = {128},
  publisher = {American Society of Civil Engineers},
}

@InProceedings{GOKTEPE-09,
  author       = {Goktepe, A. B. and Lav, A. H. and Altun, S.},
  booktitle    = {Proceedings of the Institution of Civil Engineers-Transport},
  title        = {Method for optimal vertical alignment of highways},
  year         = {2009},
  organization = {Thomas Telford Ltd},
  pages        = {177--188},
  volume       = {162},
}

@InProceedings{TAT-03,
  author       = {Tat, C. W. and Tao, F.},
  booktitle    = {Proceedings of the 2003 IEEE International Conference on Intelligent Transportation Systems},
  title        = {{Using GIS and genetic algorithm in highway alignment optimization}},
  year         = {2003},
  organization = {IEEE},
  pages        = {1563--1567},
  volume       = {2},
}

@Article{GOH-88,
  author  = {C. J. Goh and E. P. Chew and T. F. Fwa},
  journal = {Transportation Research Part B: Methodological},
  title   = {Discrete and continuous models for computation of optimal vertical highway alignment},
  year    = {1988},
  issn    = {0191-2615},
  number  = {6},
  pages   = {399-409},
  volume  = {22},
  doi     = {https://doi.org/10.1016/0191-2615(88)90021-5},
}

@Article{CHEW-89,
  author    = {Chew, E. P. and Goh, C. J. and Fwa, T. F.},
  journal   = {Transportation Research Part B: Methodological},
  title     = {Simultaneous optimization of horizontal and vertical alignments for highways},
  year      = {1989},
  number    = {5},
  pages     = {315--329},
  volume    = {23},
  publisher = {Elsevier},
}

@Article{LI-13,
  author  = {Li, W. and Pu, H. and Zhao, H. and Liu, W.},
  journal = {Journal of Software},
  title   = {{Approach for Optimizing 3D Highway Alignments Based on Two-stage Dynamic Programming.}},
  year    = {2013},
  number  = {11},
  pages   = {2967--2973},
  volume  = {8},
}

@Mastersthesis{HOSSAIN-13,
  author = {Hossain, S.},
  school = {University of British Columbia},
  title  = {Models and strategies for efficiently optimizing the vertical alignment of roads for multimaterial},
  year   = {2013},
}

@Article{GARCIA-CHAN-21,
  author  = {N. García-Chan and L. J. Alvarez-Vázquez and A. Martínez and M. E. Vázquez-Méndez},
  journal = {Mathematics and Computers in Simulation},
  title   = {Designing an ecologically optimized road corridor surrounding restricted urban areas: A mathematical methodology},
  year    = {2021},
  issn    = {0378-4754},
  pages   = {745-759},
  volume  = {190},
  doi     = {https://doi.org/10.1016/j.matcom.2021.06.016},
  url     = {https://www.sciencedirect.com/science/article/pii/S0378475421002391},
}

@Article{PUSHAK-16,
  author    = {Pushak, Y. and Hare, W. and Lucet, Y.},
  journal   = {European Journal of Operational Research},
  title     = {Multiple-path selection for new highway alignments using discrete algorithms},
  year      = {2016},
  number    = {2},
  pages     = {415--427},
  volume    = {248},
  publisher = {Elsevier},
}

@Article{HIRPA-16,
  author    = {Hirpa, D. and Hare, W. and Lucet, Y. and Pushak, Y. and Tesfamariam, S.},
  journal   = {Transportation Research Part C: Emerging Technologies},
  title     = {A bi-objective optimization framework for three-dimensional road alignment design},
  year      = {2016},
  pages     = {61--78},
  volume    = {65},
  publisher = {Elsevier},
}

@InProceedings{LOEFBERG-04,
  author    = {L{\"{o}}fberg, J.},
  booktitle = {In Proceedings of the CACSD Conference},
  title     = {{YALMIP : A Toolbox for Modeling and Optimization in MATLAB}},
  year      = {2004},
  address   = {Taipei, Taiwan},
}

@Manual{MATLAB-21,
  title     = {9.12.0.2170939 (R2022a) 64-bit (win64)},
  address   = {Natick, Massachusetts},
  author    = {MATLAB},
  year      = {2021},
  publisher = {The MathWorks Inc.},
}

@Manual{CPLEX-09,
  title   = {{V12. 1: User’s Manual for CPLEX}},
  author  = {{CPLEX, IBM ILOG}},
  year    = {2009},
  journal = {International Business Machines Corporation},
  number  = {53},
  volume  = {46},
}

@Article{MOMO-23,
  author  = {Momo, N. S. and Hare, W. and Lucet, Y.},
  journal = {Computer-Aided Civil and Infrastructure Engineering},
  title   = {Modeling side slopes in vertical alignment resource road construction using convex optimization},
  year    = {2023},
  number  = {2},
  pages   = {211-224},
  volume  = {38},
  doi     = {https://doi.org/10.1111/mice.12879},
  url     = {https://onlinelibrary.wiley.com/doi/abs/10.1111/mice.12879},
}

@Book{STARK-72,
  author    = {Stark, R. and Nicholls, R.},
  publisher = {McGraw-Hill},
  title     = {Mathematical foundations for design: Civil engineering systems},
  year      = {1972},
  address   = {New York},
}

@Article{MAYER-81,
  author  = {Mayer, R. H. and Stark, R. M.},
  journal = {Journal of the Construction Division},
  title   = {Earthmoving Logistics},
  year    = {1981},
  number  = {2},
  pages   = {297-312},
  volume  = {107},
  doi     = {10.1061/JCCEAZ.0000963},
  url     = {https://ascelibrary.org/doi/abs/10.1061/JCCEAZ.0000963},
}

@Article{AYMAN-23,
  author    = {Ayman,K. M. and Hare, W. and Lucet, Y.},
  journal   = {Engineering Optimization},
  title     = {A multi-road quasi network flow model for the vertical alignment optimization of a road network},
  year      = {2023},
  pages     = {1-24},
  doi       = {10.1080/0305215X.2023.2227841},
  publisher = {Taylor & Francis},
}

@Article{GAO-22,
  author  = {Gao, T. and Li, Z. and Gao, Y. and Schonfeld, P. and Feng, X. and Wang, Q. and He, Q.},
  journal = {Computer-Aided Civil and Infrastructure Engineering},
  title   = {A deep reinforcement learning approach to mountain railway alignment optimization},
  year    = {2022},
  number  = {1},
  pages   = {73-92},
  volume  = {37},
  doi     = {https://doi.org/10.1111/mice.12694},
}

@Book{RODRIGUE-06,
  author    = {Rodrigue, J. P. and Comtois, C. and Slack, B.},
  publisher = {Taylor \& Francis},
  title     = {The Geography of Transport Systems},
  year      = {2006},
  isbn      = {9781134257782},
  url       = {https://books.google.ca/books?id=7fuzJLvds4UC},
}

@Article{SONG-23,
  author     = {Song, T. and Schonfeld, P. and Pu, H.},
  journal    = {IEEE Transactions on Intelligent Transportation Systems},
  title      = {A Review of Alignment Optimization Research for Roads, Railways and Rail Transit Lines},
  year       = {2023},
  issn       = {1524-9050},
  number     = {5},
  pages      = {4738–4757},
  volume     = {24},
  doi        = {10.1109/TITS.2023.3235685},
  issue_date = {May 2023},
  publisher  = {IEEE Press},
  url        = {https://doi.org/10.1109/TITS.2023.3235685},
}

@InBook{ACHTERBERG-13,
  author    = {Achterberg, T. and Wunderling, R.},
  pages     = {449--481},
  publisher = {Springer Berlin Heidelberg},
  title     = {Mixed Integer Programming: Analyzing 12 Years of Progress},
  year      = {2013},
  address   = {Berlin, Heidelberg},
  isbn      = {978-3-642-38189-8},
  booktitle = {Facets of Combinatorial Optimization: Festschrift for Martin Gr{\"o}tschel},
  doi       = {10.1007/978-3-642-38189-8_18},
  url       = {https://doi.org/10.1007/978-3-642-38189-8_18},
}

@Article{GWAK-18,
  author  = {Han-Seong Gwak and Jongwon Seo and Dong-Eun Lee},
  journal = {Automation in Construction},
  title   = {Optimal cut-fill pairing and sequencing method in earthwork operation},
  year    = {2018},
  issn    = {0926-5805},
  pages   = {60-73},
  volume  = {87},
  doi     = {https://doi.org/10.1016/j.autcon.2017.12.010},
  url     = {https://www.sciencedirect.com/science/article/pii/S0926580517305009},
}

@Article{CASAL-17,
	author   = {Gerardo Casal and Duarte Santamarina and Miguel E. Vázquez-Méndez},
	journal  = {Transportation Research Part C: Emerging Technologies},
	title    = {Optimization of horizontal alignment geometry in road design and reconstruction},
	year     = {2017},
	issn     = {0968-090X},
	pages    = {261-274},
	volume   = {74},
	doi      = {https://doi.org/10.1016/j.trc.2016.11.019},
	url      = {https://www.sciencedirect.com/science/article/pii/S0968090X16302431},
}

@Article{ZHANG-23,
	author  = {Zhang, Tianlong and Gao, Yan and Gao, Tianci and Schonfeld, Paul and Wu, Yuecheng and Zhu, Ying and Yang, Shusheng and Wang, Ping and He, Qing},
	journal = {Computer-Aided Civil and Infrastructure Engineering},
	title   = {A sequential exploration algorithm for the design optimization of horizontal road alignment},
	year    = {2023},
	number  = {15},
	pages   = {2049-2071},
	volume  = {38},
	doi     = {https://doi.org/10.1111/mice.12990},
	eprint  = {https://onlinelibrary.wiley.com/doi/pdf/10.1111/mice.12990},
	url     = {https://onlinelibrary.wiley.com/doi/abs/10.1111/mice.12990},
}
	
\appendix 

\section*{Appendix A: QCQP-QNF vs. CUVA}
\label{appendix:qcqp-old-new-comp}
{
	\small
	\begin{longtable}{crrrrrr}
		\caption{Comparison of approximation errors of cut sections for QCQP-QNF and CUVA model.} \label{tab:qcqp-qnf-cut-error} \\
		\toprule
		\thead{Road \\ \#} & \thead{Sta-\\tions} & \thead{Slabs} & \thead{MAPE \\ (QCQP-\\QNF)} & \thead{MAPE \\ (CUVA)} & \thead{RMSE \\ (QCQP-\\QNF)} & \thead{RMSE \\ (CUVA)}\\
		\midrule
		\endfirsthead
		
		\multicolumn{7}{c}%
		{{\bfseries Table \thetable\ continued from previous page}} \\
		\toprule
		\thead{Road \\ \#} & \thead{Sta-\\tions} & \thead{Slabs} & \thead{MAPE \\ (QCQP-\\QNF)} & \thead{MAPE \\ (CUVA)} & \thead{RMSE \\ (QCQP-\\QNF)} & \thead{RMSE \\ (CUVA)}\\
		\midrule
		\endhead
		
		\midrule
		\multicolumn{7}{r}{{Continued on next page}} \\ 
		\endfoot
		
		\endlastfoot
		
		1       & 10       & 22    & 8.41         & \textbf{3.62}      & 6.37          & \textbf{2.80}             \\
		2       & 14       & 22    & 8.23          & \textbf{3.54}           & 5.89         & \textbf{2.57}               \\
		3       & 18       & 22    & 8.32          & \textbf{3.59}            & 6.20         & \textbf{2.71}   \\
		4       & 26       & 22    & 7.44         & \textbf{3.34}   & 5.92         & \textbf{2.70}   \\
		5       & 28       & 22    & 7.31           & \textbf{3.30}  & 5.93         & \textbf{2.71}               \\
		6       & 45       & 22    & 8.84          & \textbf{3.64}   & 7.66         & \textbf{3.01}            \\
		7       & 45       & 22    & 8.87          & \textbf{3.58}   & 7.73          & \textbf{3.03}           \\
		8       & 45       & 22    & 7.41          & \textbf{3.30}   & 5.57         & \textbf{2.54}      \\
		9       & 48       & 22    & 7.62          & \textbf{3.28}   & 6.43         & \textbf{2.69}        \\
		10      & 50       & 22    & 6.96          & \textbf{3.18}   & 5.56         & \textbf{2.60}      \\
		11      & 50       & 22    & 8.20          & \textbf{3.44}   & 6.55         & \textbf{2.71}        \\
		12      & 51       & 29    & 0.01          & \textbf{0.00}   & 0.00          & \textbf{0.00}      \\
		13      & 51       & 11    & 0.59          & \textbf{0.48}   & 0.52          & \textbf{0.43}        \\
		14      & 51       & 27    & 8.67          & \textbf{3.24}   & 7.09          & \textbf{2.63}        \\
		15      & 52       & 27    & 7.92          & \textbf{3.10}   & 5.80         & \textbf{2.31}         \\
		16      & 53       & 27    & 7.48          & \textbf{3.01}   & 5.60         & \textbf{2.29}      \\
		17      & 54       & 22    & 7.43          & \textbf{3.34}   & 5.75         & \textbf{2.63}     \\
		18      & 54       & 22    & 7.49          & \textbf{3.34}   & 5.75          & \textbf{2.60}     \\
		19      & 57       & 22    & 8.05          & \textbf{3.45}   & 5.46         & \textbf{2.40}       \\
		20      & 60       & 22    & 7.56          & \textbf{3.37}   & 5.90         & \textbf{2.68}        \\
		21      & 61       & 22    & 7.68         & \textbf{3.39}    & 5.82         & \textbf{2.61}        \\
		22      & 63       & 22    & 7.63          & \textbf{3.41}   & 6.14          & \textbf{2.79}       \\
		23      & 68       & 22    & 7.62          & \textbf{3.35}   & 5.46          & \textbf{2.46}        \\
		24      & 69       & 27    & 8.31          & \textbf{3.17}   & 5.96         & \textbf{2.31}         \\
		25      & 71       & 36    & 7.97           & \textbf{2.65}  & 5.62          & \textbf{1.89}        \\
		26      & 82       & 22    & 7.50         & \textbf{3.34}   & 5.79          & \textbf{2.64}         \\
		27      & 85       & 22    & 7.75          & \textbf{3.41}   & 5.80          & \textbf{2.61}       \\
		28      & 89       & 27    & 7.91          & \textbf{3.11}    & 5.91          & \textbf{2.36}       \\
		29      & 95       & 22    & 8.19           & \textbf{3.42}   & 6.26            & \textbf{2.57}     \\
		30      & 101      & 41    & 0.02          & \textbf{0.01}    & 0.00          & \textbf{0.00}        \\
		31      & 101      & 25    & 0.02          & \textbf{0.01}    & 0.00          & \textbf{0.00}        \\
		32      & 132      & 22    & 8.22          & \textbf{3.50}    & 7.00           & \textbf{2.95}       \\
		33      & 151      & 23    & 0.01          & \textbf{0.00}    & 0.00          & \textbf{0.00}        \\
		34      & 151      & 25    & 0.04          & \textbf{0.02}    & 0.01          & \textbf{0.00}         \\
		35      & 201      & 17    & 0.02          & \textbf{0.01}    & 0.01          & \textbf{0.00}        \\
		36      & 202      & 101   & 0.02          & \textbf{0.00}    & 0.00          & \textbf{0.00}        \\
		37      & 451      & 29    & 0.05          & \textbf{0.02}    & 0.02            & \textbf{0.01}    \\
		38      & 971      & 30    & 8.53          & \textbf{3.01}    & 7.16          & \textbf{2.48}       \\
		39      & 1,001     & 21    & 0.02          & \textbf{0.01}   & 0.00          & \textbf{0.00}       \\
		40      & 1,056     & 30    & 7.86          & \textbf{2.88}   & 5.79          & \textbf{2.15}        \\
		41      & 1,092     & 31    & 7.86          & \textbf{2.89}   & 5.81          & \textbf{2.16}       \\
		42      & 1,281     & 31    & 8.02          & \textbf{2.92}   & 6.16          & \textbf{2.25}        \\
		43      & 1,290     & 31    & 7.94           & \textbf{2.86}  & 5.84         & \textbf{2.10}        \\
		44      & 1,506     & 31    & 7.56          & \textbf{2.79}   & 5.30          & \textbf{1.98}       \\
		45      & 1,732     & 31    & 7.68          & \textbf{2.84}   & 5.53          & \textbf{2.07}      \\
		46      & 2,286     & 11    & 0.01          & \textbf{0.01}   & 0.00          & \textbf{0.00}     \\
		47      & 2,118     & 30    & 7.81          & \textbf{2.88}   & 5.85          & \textbf{2.18}       \\
		48      & 2,247     & 22    & 16.21         & \textbf{5.82}   & 163.06        & \textbf{51.94}      \\
		49      & 3,704     & 31    & 7.99          & \textbf{2.92}   & 5.94         & \textbf{2.20}      \\
		50      & 4,194     & 30    & 8.32         & \textbf{2.88}   & 6.71         & \textbf{2.24}       \\
		51      & 4,325     & 35    & 12.75         & \textbf{3.01}  & 16.83         & \textbf{3.31}     \\
		52      & 13,892    & 12    & 6.94          & \textbf{2.29}  & 5.18         & \textbf{1.58}        \\ \hline
		\bottomrule
	\end{longtable}
}

{
	\small
	\begin{longtable}{crrrrrr}
		\caption{Comparison of approximation errors of fill sections for QCQP-QNF and CUVA model.} \label{tab:qcqp-qnf-fill-error} \\
		\toprule
		\thead{Road \\ \#} & \thead{Sta-\\tions} & \thead{Slabs} & \thead{MAPE \\ (QCQP-\\QNF)} & \thead{MAPE \\ (CUVA)} & \thead{RMSE \\ (QCQP-\\QNF)} & \thead{RMSE \\ (CUVA)}\\
		\midrule
		\endfirsthead
		
		\multicolumn{7}{c}%
		{{\bfseries Table \thetable\ continued from previous page}} \\
		\toprule
		\thead{Road \\ \#} & \thead{Sta-\\tions} & \thead{Slabs} & \thead{MAPE \\ (QCQP-\\QNF)} & \thead{MAPE \\ (CUVA)} & \thead{RMSE \\ (QCQP-\\QNF)} & \thead{RMSE \\ (CUVA)}\\
		\midrule
		\endhead
		
		\midrule
		\multicolumn{7}{r}{{Continued on next page}} \\ 
		\endfoot
		
		\endlastfoot
		
		1       & 10       & 22    & 2.14          & \textbf{1.56}               & 2.42         & \textbf{1.75}              \\
		2       & 14       & 22    & 1.37           & \textbf{1.11}         & 2.15         & \textbf{1.42}             \\
		3       & 18       & 22    & 1.85          & \textbf{1.60}        & 3.08          & \textbf{2.66}               \\
		4       & 26       & 22    & 1.97          & \textbf{1.56}          & 2.62          & \textbf{2.17}         \\
		5       & 28       & 22    & 2.29          & \textbf{1.80}      & 2.78          & \textbf{2.55}          \\
		6       & 45       & 22    & 2.40          & \textbf{2.16}      & 4.70          & 4.52    \\
		7       & 45       & 22    & 2.36          & \textbf{1.89}           & 4.58          & \textbf{4.27}            \\
		8       & 45       & 22    & 1.84          & \textbf{1.52}        & 2.35          & \textbf{2.17}           \\
		9       & 48       & 22    & 2.25          & \textbf{2.03}            & \textbf{3.08}         & 3.36            \\
		10      & 50       & 22    & 2.00          & \textbf{1.26}           & 1.88          & \textbf{1.18}           \\
		11      & 50       & 22    & 2.31          & \textbf{1.57}             & 2.76          & \textbf{2.18}           \\
		12      & 51       & 29    & 0.01          & \textbf{0.00}            & 0.00          & \textbf{0.00}          \\
		13      & 51       & 11    & 0.26          & \textbf{0.22}         & 0.19          & \textbf{0.16}             \\
		14      & 51       & 27    & 1.97          & \textbf{1.54}        & 2.56          & \textbf{2.19}         \\
		15      & 52       & 27    & 2.17          & \textbf{1.35}       & 2.31          & \textbf{1.52}            \\
		16      & 53       & 27    & 2.29          & \textbf{1.37}      & 1.99         & \textbf{1.28}         \\
		17      & 54       & 22    & 1.72          & \textbf{1.20}        & 2.27          & \textbf{1.50}            \\
		18      & 54       & 22    & 2.12          & \textbf{1.33}        & 2.14          & \textbf{1.39}         \\
		19      & 57       & 22    & 3.07          & \textbf{1.87}          & 2.35          & \textbf{1.68}          \\
		20      & 60       & 22    & 2.37          & \textbf{1.80}         & 2.63          & \textbf{1.98}            \\
		21      & 61       & 22    & 2.25          & \textbf{1.48}        & 2.24          & \textbf{1.49}        \\
		22      & 63       & 22    & 2.25          & \textbf{1.43}      & 2.23          & \textbf{1.47}        \\
		23      & 68       & 22    & 2.41          & \textbf{1.47}         & 2.05          & \textbf{1.29}           \\
		24      & 69       & 27    & 2.37          & \textbf{1.45}       & 2.27          & \textbf{1.44}          \\
		25      & 71       & 36    & 2.50          & \textbf{1.27}       & 1.99          & \textbf{1.07}             \\
		26      & 82       & 22    & 2.00          & \textbf{1.41}          & 1.83          & \textbf{1.44}           \\
		27      & 85       & 22    & 2.41          & \textbf{1.62}         & 2.39          & \textbf{1.74}             \\
		28      & 89       & 27    & 2.41          & \textbf{1.42}      & 2.08         & \textbf{1.33}          \\
		29      & 95       & 22    & 2.57          & \textbf{1.68}        & 2.13          & \textbf{1.60}            \\
		30      & 101      & 41    & 0.02          & \textbf{0.01}        & 0.00          & \textbf{0.00}           \\
		31      & 101      & 25    & 0.02          & \textbf{0.01}          & 0.00          & \textbf{0.00}           \\
		32      & 132      & 22    & \textbf{73.77}        & 169.82          & 7.24          & \textbf{7.11}             \\
		33      & 151      & 23    & 0.01          & \textbf{0.00}           & 0.00          & \textbf{0.00}           \\
		34      & 151      & 25    & 0.04         & \textbf{0.02}        & 0.01          & \textbf{0.00}            \\
		35      & 201      & 17    & 0.02          & \textbf{0.01}      & 0.01          & \textbf{0.00}             \\
		36      & 202      & 101   & 0.02          & \textbf{0.00}        & 0.00          & \textbf{0.00}            \\
		37      & 451      & 29    & 0.06          & \textbf{0.03}       & 0.01          & \textbf{0.01}           \\
		38      & 971      & 30    & 2.00          & \textbf{1.57}      & 2.65         & \textbf{2.43}             \\
		39      & 1,001     & 21    & 0.02          & \textbf{0.01}     & 0.00          & \textbf{0.00}             \\
		40      & 1,056     & 30    & 2.18          & \textbf{1.53}          & 2.50          & \textbf{1.76}            \\
		41      & 1,092     & 31    & 2.45          & \textbf{1.43}       & 2.18          & \textbf{1.31}             \\
		42      & 1,281     & 31    & 2.24          & \textbf{1.36}       & 2.23          & \textbf{1.45}           \\
		43      & 1,290     & 31    & 2.68         & \textbf{1.59}         & 2.38          & \textbf{1.56}       \\
		44      & 1,506     & 31    & 2.46          & \textbf{1.36}        & 2.00          & \textbf{1.18}       \\
		45      & 1,732     & 31    & 2.16          & \textbf{1.45}       & 2.27          & \textbf{1.76}         \\
		46      & 2,286     & 11    & 0.01          & \textbf{0.00}     & 0.00          & \textbf{0.00}          \\
		47      & 2,118     & 30    & 2.43          & \textbf{1.38}        & 2.09         & \textbf{1.27}       \\
		48      & 2,247     & 22    & \textbf{6.96}          & 11.90        & 49.57         & \textbf{46.72}      \\
		49      & 3,704     & 31    & 2.59          & \textbf{1.65}          & 2.58         & \textbf{1.76}         \\
		50      & 4,194     & 30    & 2.44          & \textbf{2.09}       & \textbf{3.20}          & 3.42             \\
		51      & 4,325     & 35    & 3.09          & \textbf{2.49}        & \textbf{17.4}        & 20.57           \\
		52      & 13,892    & 12    & 2.16          & \textbf{1.33}        & 2.20         & 2.15         \\ \hline
		\bottomrule
	\end{longtable}
}

\section*{Appendix B: MH-QNF vs. CUVA}
\label{appendix:qcqp-milp-comp}
{
	\small
	\begin{longtable}{cccrrrrrr}
		\caption{Comparison of optimal costs for MH-QNF and CUVA model.} \label{tab:MH-QNF-optimal-cost-comparison} \\
		\toprule
		\thead{Road \\ \#} & \thead{Sta-\\tions} & \thead{Slabs} & \thead{Optimal \\Cost \\ (MH-\\QNF)} & \thead{Optimal \\Cost \\ (CUVA)} & \thead{Diff \\ (\%)} & \thead{Seconds \\ (MH-\\QNF)} & \thead{Seconds \\ (CUVA)} & \thead{Speed-\\up}  \\
		\midrule
		\endfirsthead
		
		\multicolumn{9}{c}%
		{\tablename\ \thetable\ -- \textit{Continued from previous page}} \\
		\toprule
		\thead{Road \\ \#} & \thead{Sta-\\tions} & \thead{Slabs} & \thead{Optimal \\Cost \\ (MH-\\QNF)} & \thead{Optimal \\Cost \\ (CUVA)} & \thead{Diff \\ (\%)} & \thead{Seconds \\ (MH-\\QNF)} & \thead{Seconds \\ (CUVA)} & \thead{Speed-\\up}  \\
		\midrule
		\endhead
		
		\midrule
		\multicolumn{9}{r}{\textit{Continued on next page}} \\
		\endfoot
		
		\bottomrule
		\endlastfoot
		1       & 51       & 10    & 217                 & \textbf{146}                     & \textbf{33} & \textbf{0.03}             & 0.07                    & 0.43 \\
		2       & 51       & 28    & 197               & \textbf{148}                     & \textbf{25} & \textbf{0.03}             & 0.07                    & 0.42 \\
		3       & 51       & 30    & 197               & \textbf{147}                      & \textbf{25} & \textbf{0.03}             & 0.07                    & 0.42 \\
		4       & 51       & 40    & 191               & \textbf{148}                     & \textbf{23} & \textbf{0.03}             & 0.07                    & 0.42 \\
		5       & 51       & 10    & 567               & \textbf{463}                     & \textbf{18} & \textbf{0.02}             & 0.08                    & 0.25     \\
		6       & 51       & 40    & 509               & \textbf{465}                     & \textbf{9} & \textbf{0.03}             & 0.08                    & 0.37    \\
		7       & 101      & 10    & Infeasible            & \textbf{150}                     & -          & Infeasible       & \textbf{0.10}                     & -       \\
		8       & 101      & 17    & 364               & \textbf{151}                     & \textbf{58} & \textbf{0.04}             & 0.11                    & 0.36 \\
		9       & 101      & 24    & 353               & \textbf{153}                     & \textbf{57} & \textbf{0.05}             & 0.1                     & 0.50      \\
		10      & 151      & 10    & 694               & \textbf{499}                     & \textbf{28} & \textbf{0.05}             & 0.12                    & 0.41 \\
		11      & 151      & 6     & 731               & \textbf{498}                     & \textbf{32} & \textbf{0.05}             & 0.11                    & 0.45 \\
		12      & 151      & 16    & 663               & \textbf{500}                     & \textbf{25}  & \textbf{0.06}             & 0.11                    & 0.54 \\
		13      & 151      & 22    & 650               & \textbf{501}                     & \textbf{23} & \textbf{0.06}             & 0.11                    & 0.54 \\
		14      & 151      & 10    & 3,130               & \textbf{2,820}                       & \textbf{10} & \textbf{0.05}             & 0.11                    & 0.45 \\
		15      & 151      & 24    & 2,951               & 2,820                       & 4 & \textbf{0.06}             & 0.12                    & 0.50      \\
		16      & 201      & 10    & 3,808                & \textbf{3,427}                         & \textbf{10} & \textbf{0.06}             & 0.12                    & 0.50      \\
		17      & 201      & 6     & 4,042               & \textbf{3,427}                       & \textbf{15} & \textbf{0.05}             & 0.13                    & 0.38 \\
		18      & 201      & 16    & 3,650               & \textbf{3,427}                       & \textbf{6} & \textbf{0.07}             & 0.14                    & 0.50      \\
		19      & 451      & 10    & 5,018               & \textbf{4,210}                       & \textbf{16} & \textbf{0.15}             & 0.21                    & 0.71 \\
		20      & 451      & 4     & 5,503               & \textbf{4,209}                       & \textbf{24} & \textbf{0.12}             & 0.22                    & 0.54 \\
		21      & 451      & 16    & 4,812               & \textbf{4,211}                       & \textbf{12} & \textbf{0.17}             & 0.28                    & 0.60 \\
		22      & 451      & 20    & 4,752               & \textbf{3,958}                       & \textbf{17} & \textbf{0.21}             & 0.24                    & 0.87    \\
		23      & 451      & 30    & 4,693              & \textbf{4,217}                       & \textbf{10} & \textbf{0.22}             & 0.29                    & 0.75 \\ \hline
	\end{longtable}

}
	
\end{document}